\documentclass[a4paper,12pt]{article}
\usepackage{graphicx,epsfig}
\usepackage[T1,T2A]{fontenc}
\usepackage[utf8]{inputenc}
\usepackage{amsfonts, amsthm}
\usepackage{xcolor}

\newtheorem{proposition}{Proposition}

\newcommand{\ptl}{\partial}

\title{Asymptotic evaluation of three-dimensional integrals with singularities in application to wave phenomena}
\author{А.~V. Shanin, A.~Yu. Laptev}

\begin{document}

\maketitle

\begin{abstract}
We consider a three-dimensional Fourier integral in which the exponent in the exponential factor is the product of some phase function and a large parameter. The asymptotics of this integral is sought when the large parameter tends to infinity. In the one-dimensional case, the asymptotics of such an integral is constructed by the points of stationary phase and singularities of the integrand. The three-dimensional case is more complicated: special points such as points of stationary phase in the domain, on singularity, on the crossing of singularities, points of triple crossing of singularities, and also conical points of the singularities, can contribute to the asymptotics. For all these types of singularities, topological conditions for the existence of nonzero asymptotics are constructed, and the asymptotics themselves are derived. The proposed technique is tested on the example of the classical problem of Kelvin waves on the surface of a deep fluid behind a towed body.
\end{abstract}

\section{Introduction}

Typically, wave phenomena can be described by integrals close in form to multidimensional Fourier integrals. In particular, such integrals can be plane wave decompositions (spectral representations). Accordingly, understanding of wave phenomena requires a technique for constructing asymptotic estimations of such integrals. For one-dimensional integrals, one can use the stationary phase method proposed by Kelvin (see the discussion in~\cite{Martin2025}) or the related method of steepest descent. For multivariate integrals, the multivariate method of steepest descent described in~\cite{Fedoryuk1977} applies in some cases, but there are problems to which the standard methods cannot be applied. These include integrals whose integrands contain singularities. For such problems in the two-dimensional case, a method~\cite{Mironov2021,Assier2022,Shanin2024} has been proposed, generalizing the stationary phase method in Poincaré's form~\cite{Poincare1904} to two-dimensional integrals. The purpose of this paper is to generalize this method to the three-dimensional case. The idea of the method remains the same, but the number of types of special points increases significantly.

As is known, the Poincaré's work proposes to shift slightly the integration contour from the real axis into the complex domain so that the integrand becomes exponentially small almost everywhere. The stationary phase point is thus transformed into a steepest descent point (below we do not distinguish between stationary phase points and steepest descent points), and singularities become bypassed by small loops. We perform a similar procedure in the three-dimensional case. Naturally, in complex 3D, topological difficulties arise in describing the deformation of the integration domain and in checking that the integration domain does not cross the singularity. This is the main focus of this paper.

The advantage of our approach is that, although the deformations of the integration domain are quite complicated and occur in three-dimensional {\em complex\/} (i.e., six-dimensional) space, the special points can be analyzed by considering only the intersections of singularities with three-dimensional {\em real\/} space, i.e., quite easily.

We show that in the general case (for almost all values of the parameters) the principle of locality works: the asymptotic estimate of the integral is constructed as the sum of several contributions from special points. Our consideration of two-dimensional problems, as well as standard one-dimensional analysis, shows that this kind of estimation has not only a mathematical but also a physical meaning. Namely, each of the contributions obtained from special points (up to some obvious symmetries) corresponds to 
one of the types of the waves supported by the system.

The paper is motivated by the classical problem of Kelvin waves on the surface of a deep liquid. For the first time, this problem appears in the lecture given by Kelvin in 1887 (the published version is in \cite{Thomson1887}). Although there are no formulas in the published version of this lecture, it is stated that the wave pattern formed in the rear of the ship is located inside two lines drawn from the bow of the ship and inclined to the wake on its two sides at equal angles of $\arctan{(1/2\sqrt{2})\approx19^\circ 28'}$. This angle is known as the Kelvin angle.

Later, in 1906, Kelvin wrote a paper \cite{Kelvin1906} in which he analyzed the ship waves in more detail. In particular, the first integral representation of the displacement of the water surface by a disturbance traveling with constant velocity was given. This study was later continued in~\cite{Lamb1916a,Lamb1916b}.

In book~\cite{Stoker1957}, the transient water waves are investigated in detail. The formulas for velocity potential and surface displacement caused by disturbances that originated at the water surface are given for the case of infinite depth of water. The formulas are analyzed by the method of the stationary phase. Linearization of the problem is also provided. Similar analysis is provided in~\cite{Wehausen1960}.

In paper~\cite{Liu2001}, the transient ship waves induced by a moving body are analyzed. Using the method proposed in \cite{Wehausen1960}, the integral representation of velocity potential and surface displacement in a coordinate system moving with the ship is given. Both potential and displacement are represented as the sum of two terms: the first is the steady-state term, the second is the transient term.

The current paper is organized as follows. In Section~2, the types of integrals and the basic assumptions about these integrals, which should be satisfied for our approach to work, are formulated. In Section~3, the concepts of {\em desired} and {\em allowed} deformations of real three-dimensional integration domains in a complex space are introduced. In Section~4, the 
indentation of the integration domain from the singularities is described. In Section~5, the principle of locality (i.e., the principle by which the full integral can be reduced to the sum of contributions from several special points) is formulated. In Section~6, the points that {\em do not\/} contribute to the integral estimate are described. In Section~7, the points that {\em do\/} contribute to the integral estimate are described. The leading terms of the asymptotics of these contributions are derived. In Section~8, the contribution from the boundary of the integration domain is discussed. In Section~9, the approach described above is used to find the asymptotics of the three-dimensional integral representation of displacement caused by the body moving on the surface of fluid with constant velocity.


\section{Problem statement and basic assumptions
}

Consider the integral 
\begin{equation}
u(z; \Lambda) = \int_\Gamma F(\xi) \, \exp \{ i \Lambda G(\xi ; z)\} \, d\xi,  
\label{e:001}
\end{equation}
where $\xi = (\xi_1, \xi_2, \xi_3)$ is the triplet of complex integration variables, $z = (z_1, z_2, z_3)$ are spatial or temporal real coordinates playing the role of parameters in the integral, 
$\Lambda$ is a large real positive parameter, $\Gamma$ is a sufficiently regular oriented domain of real dimension~3 (the integration domain), 
\begin{equation}
d\xi = d\xi_1 \wedge d\xi_2 \wedge d\xi_3.
\label{e:002}
\end{equation}
Our task is to evaluate the integral (\ref{e:001}) as $\Lambda \to \infty$, retaining only  terms that do not exponentially decrease with $\Lambda$. The function $G$ will be called the {\em phase function\/} of the integral.
The asymptotic consideration will be carried out with fixed parameters $z$, so for brevity, we will replace $u(z ; \Lambda)$ by $u(\Lambda)$ and $G(\xi ; z)$ by~$G(\xi)$.

In the simplest case, the integral is initially taken over the real three-dimensional domain, which gives
$\Gamma = \mathbb{R}^3$, $d\xi = d\xi_1 \, d\xi_2 \, d\xi_3$. We will consider the deformations of $\Gamma$ in the complex domain, which forces us to use the form (\ref{e:002}). As it is known~\cite{Shabat2},
this form is the volume of a parallelepiped stretched over a triplet of small complex vectors $(v^1, v^2, v^3)$ in the space~$\xi$:
\begin{equation}
d\xi(v^1, v^2, v^3) = dV = {\rm det}\left( \begin{array}{ccc}
v_1^1 & v_2^1 & v_3^1 \\
v_1^2 & v_2^2 & v_3^2 \\
v_1^3 & v_2^3 & v_3^3 
\end{array} \right),
\label{e:003}
\end{equation}
where the elements of the matrix are the corresponding complex coordinates of the vectors. 
This formula is important not only from the theoretical point of view but also allows to perform numerical integration in complex space.

Introduce a small neighborhood of the real space in the space $\mathbb{C}^3$ of variables $\xi$: 
\[
\mathbb{R}^3_\delta = \{ \xi \in \mathbb{C}^3 \, : \, | {\rm Im} [\xi_{1,2,3}] | < \delta\}. 
\]
for some small positive $\delta$. Let us assume that $\Gamma \subset \mathbb{R}^3_\delta$.

Let us make a number of assumptions about the integral (\ref{e:001}). We will assume that:

\begin{itemize}
\item 
The functions  $F(\xi)$ and $G(\xi)$ are analytic in $\mathbb{R}^3_\delta$, i.e.\ they are holomorphic in this set except for singularities (pole sets and
branch sets) and can be continued analytically about the branch sets. 
In the examples, we will take  algebraic functions of~$\xi$ as $F$ and~$G$ .

\item
The singularities of the functions $F(\xi)$ and $G(\xi)$ lie on {\em analytic sets\/} of complex dimension~2 and real dimension~4. In other words, the irreducible components of the singularities, which we will denote as $\sigma_j$, are defined by equations 
\begin{equation}
\sigma_j = \{ \xi \in \mathbb{C}^3 \, : \, g_j (\xi) = 0 \},
\label{e:004}
\end{equation}
where $g_j$ are functions holomorphic in $\mathbb{R}^3_\delta$.
Note that only the singularities of the function~$F$ will be of practical interest. The complex equation $g_j (\xi) = 0$  represents two real equations (on the real and the imaginary parts), which gives real codimension~2 for the singularity.

\item 
Singularities have the {\em real property\/}, namely, functions $g_j$ introduced above are real for real~$\xi$. As a consequence, the intersection of the components of the singularity set 
$\sigma_j$ with real space $\mathbb{R}^3$
\[
\sigma'_j \equiv \sigma_j \cap \mathbb{R}^3 
\]
are sets of dimension~2, i.e.\ surfaces.

\item
Almost everywhere we will assume that the singularities $\sigma_j$ are 
regular sets in $\mathbb{R}^3_\delta$. This is guaranteed by
condition 
\begin{equation}
\nabla g_j(\xi) \ne {\bf 0}, \qquad \xi \in \mathbb{R}^3_\delta.
\label{e:005a}
\end{equation}
The gradient is everywhere defined as
\[
\nabla f \equiv \left( \frac{\ptl f}{\ptl \xi_1}, 
\frac{\ptl f}{\ptl \xi_2},
\frac{\ptl f}{\ptl \xi_3}
\right).
\]
We remind that for a holomorphic function $f(\xi)$, the same values are obtained if we consider partial derivatives in $\mathbb{C}^3$ by complex variables, or if we restrict the function to $\mathbb{R}^3$ and consider derivatives by real variables.

An exception to rule (\ref{e:005a}) will be a consideration of a conical singular point.

\item 
The initial integration domain $\Gamma$ is close to the real domain
$\mathbb{R}^3$. Moreover, $\Gamma$ does not intersect with the singularities~$\sigma_j$. 
 
\item 
the function $G(\xi )$ is real or has a positive imaginary part at the points $\xi \in \mathbb{R}^3$. 

\item 
The functions $F$ and $e^{i \Lambda G}$ decay at infinity fast enough, so we do not worry about the convergence of the integral.

\end{itemize}

The second condition is a general case of the behavior of singularities of holomorphic functions, and the third condition imposes an essential restriction. In particular, from the third condition it follows that $\nabla g_j(\xi)$ is a real vector for real~$\xi$. 

The assumption about the analyticity of the functions included in the integrand allows us to apply the multivariate Cauchy theorem. That is, we will use the fact that if a continuous deformation of the integration domain does not hit the singularities $\sigma_j$, then the value of the integral over the deformed domain is the same as over the original domain. Integration domains here play the same role as integration contours in one-dimensional complex analysis. 


As an example of the function $G(\xi)$ we can consider
\begin{equation}
G(\xi) = G(\xi ; z) = \xi_1 z_1 + \xi_2 z_2 + \xi_3 z_3 . 
\label{e:006}
\end{equation}
This choice corresponds to the usual three-dimensional Fourier integral.

As examples of the function $F(\xi)$ we can consider 
\[
F(\xi) = \frac{1}{\xi_1 - \xi_2 - \xi_3}
\]
for the polar set and 
\[ 
F(\xi) = \frac{1}{\sqrt{\xi_1^2  - \xi_2^2  - \xi^3_2}}
\]
for the branch set. In the first example, there is a singulatity $\sigma$ defined by the function
\[
g(\xi) = \xi_1 - \xi_2 - \xi_3,
\]
and in the second one: 
\[
g(\xi) = \xi_1^2 - \xi_2^2 - \xi_3^2.
\]
Obviously, in both cases $\sigma$ are four-dimensional manifolds and the intersections of $\sigma \cap \mathbb{R}^3$ are two-dimensional surfaces.


\section{Scheme of the method. Description of the deformation of the integration domain}

To obtain the asymptotic evaluation of the integral (\ref{e:001}) we will follow the papers \cite{Mironov2021,Assier2022,Shanin2024}. Let us construct a continuous deformation of the domain $\Gamma$ into some new domain~$\Gamma' \subset \mathbb{R}^3_\delta$ such that
\begin{equation}
{\rm Im}[G(\xi)] > 0
\label{e:008}
\end{equation}
at all $\xi \in \Gamma'$ except in the neighborhoods of a few {\em special points\/} (these points will be determined later). 
Indeed, we deform the domain in such a way that the singularities of the integrand function are not hit in the process of deformation. Obviously, only the special points' neighborhoods contribute to the evaluation of the integral for large $\Lambda$.

Let us describe the integration domain $\Gamma$ as a set of points parameterized by the real parts of $\xi_j$, namely
\begin{equation}
\Gamma = \{ \xi \in \mathbb{R}^3_\delta \, : \, {\rm Im}[\xi_j] = \eta_j ({\rm Re}[\xi_1],{\rm Re}[\xi_2],{\rm Re}[\xi_3]) , \, \, j = 1,2,3  \}.
\label{e:009}
\end{equation}
where $\eta_j$ are small continuous functions of three real variables. The functions $\eta_j$ describe the shift of the integration domain from~$\mathbb{R}^3$. 
Obviously, it makes sense to deform only the regions of $\Gamma$ where the initial values of
${\rm Im}[G]$ are not positive.

For convenience, we assume that the functions $\eta_j$ are close to being locally constant, i.e., that the derivatives of $\eta_j$ by the coordinates are small compared to the functions themselves. For details, we refer to \cite{Assier2022}, where the notion of a {\em flattable\/} integration surface is correctly introduced.

We  combine the functions $\eta_j$ into a vector: 
\[
\eta({\rm Re}[\xi]) = (\eta_1 , \eta_2, \eta_3).
\]

To describe the deformation of the integration domain, we introduce a real variable $\tau \in [0,1]$ and consider functions 
$\eta$ continuously depending on $\tau$ as a parameter: 
\[
\eta_j = \eta_j ({\rm Re}[\xi] ; \tau).
\]
This makes it possible to define a continuous family of integration domains $\Gamma(\tau)$. 
Let us assume that at the initial moment of deformation $\Gamma(0)  = \Gamma$, and at the final moment $\Gamma(1) = \Gamma'$.

The deformations that lead to the fulfillment of inequality (\ref{e:008}), will be called {\em desired}, 
and the deformations during which the integration domain does not intersect singularities of integrand functions will be called {\em allowed}. Special points are the points in the neighborhood of which there exists no desired and allowed deformation of~$\Gamma$.

If at a real point $\xi$ the value of  $G(\xi)$ is real and the gradient $\nabla G(\xi)$ is non-zero, then condition (\ref{e:008}) for small $\eta$ is rewritten as 
\begin{equation}
\nabla G(\xi) \cdot \eta(\xi) > 0
\label{e:008a}
\end{equation}
(the dot stands for the scalar product of vectors, calculated by the usual rule: ${\bf a} \cdot {\bf b} = a_1 b_1 + a_2 b_2 + a_3 b_3$).


\section{Bypass of the singularities of $\sigma_j$ by the initial domain $\Gamma$}

Let the point $\xi^* \in \mathbb{R}^3$ (the asterisk {\em does not\/} denote a complex conjugation) belong to the singularity component $\sigma_j$ and not to any other component. Let $D$ be a neighborhood of $\xi^*$ in $\mathbb{C}^3$. Consider how a domain $\Gamma$ can be located near $\sigma_j$ in~$D$. First, let us show that we can construct a domain that does not intersect~$\sigma_j$. 

Let be $\ptl g_j / \ptl \xi_1 \ne 0$ in $D$. Substitute the variables in $D$ 
\[
(\xi_1 , \xi_2, \xi_3)  \to (w_1, w_2, w_3), 
\]
where
\[
w_1 = g_j (\xi_1, \xi_2, \xi_3), \qquad w_2 = \xi_2 - \xi^*_2, \qquad w_3 = \xi_3 - \xi_3^*.
\]
In the new coordinates, the real domain is mapped to the real domain and the singularity $\sigma_j$ takes the form $w_1 = 0$.

\begin{figure}[ht]
\centerline{\epsfig{width = 10 cm, file=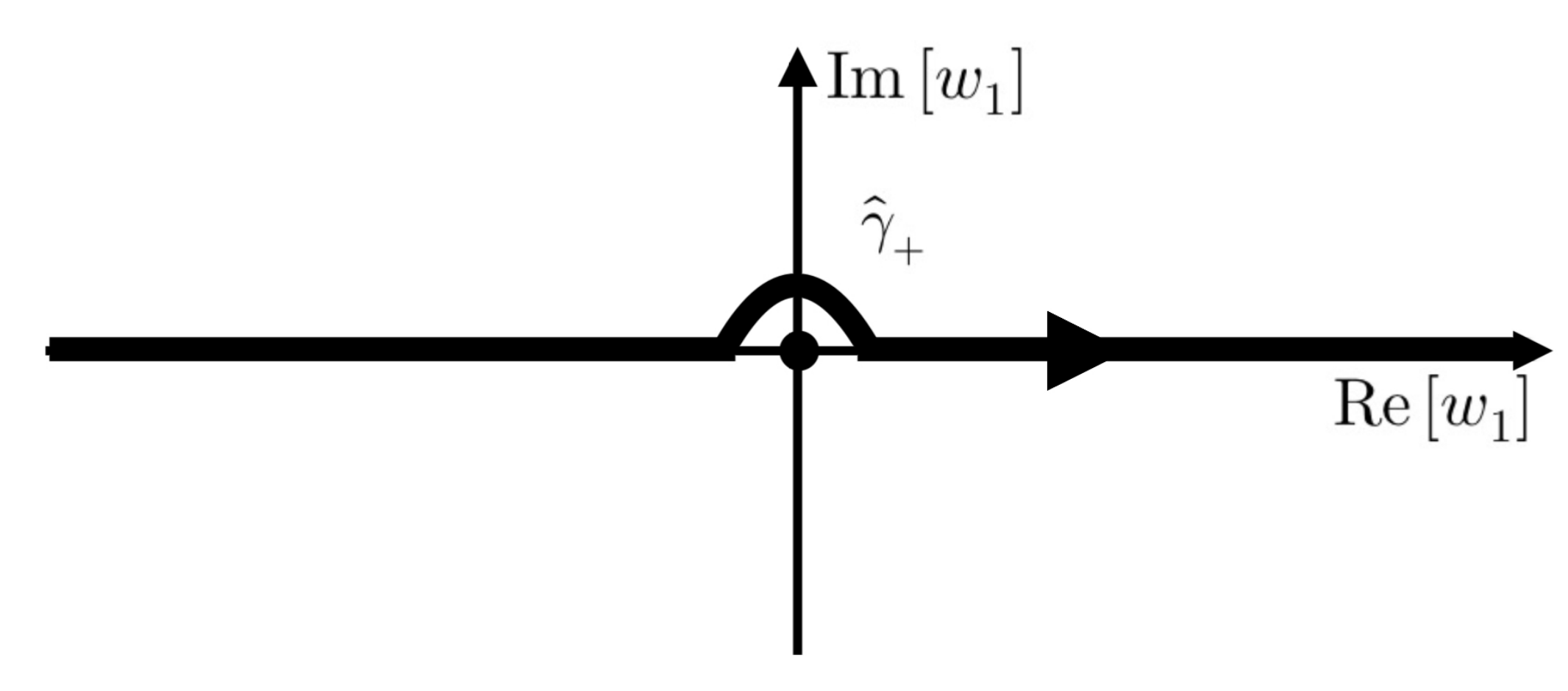}}
\centerline{\epsfig{width = 10 cm, file=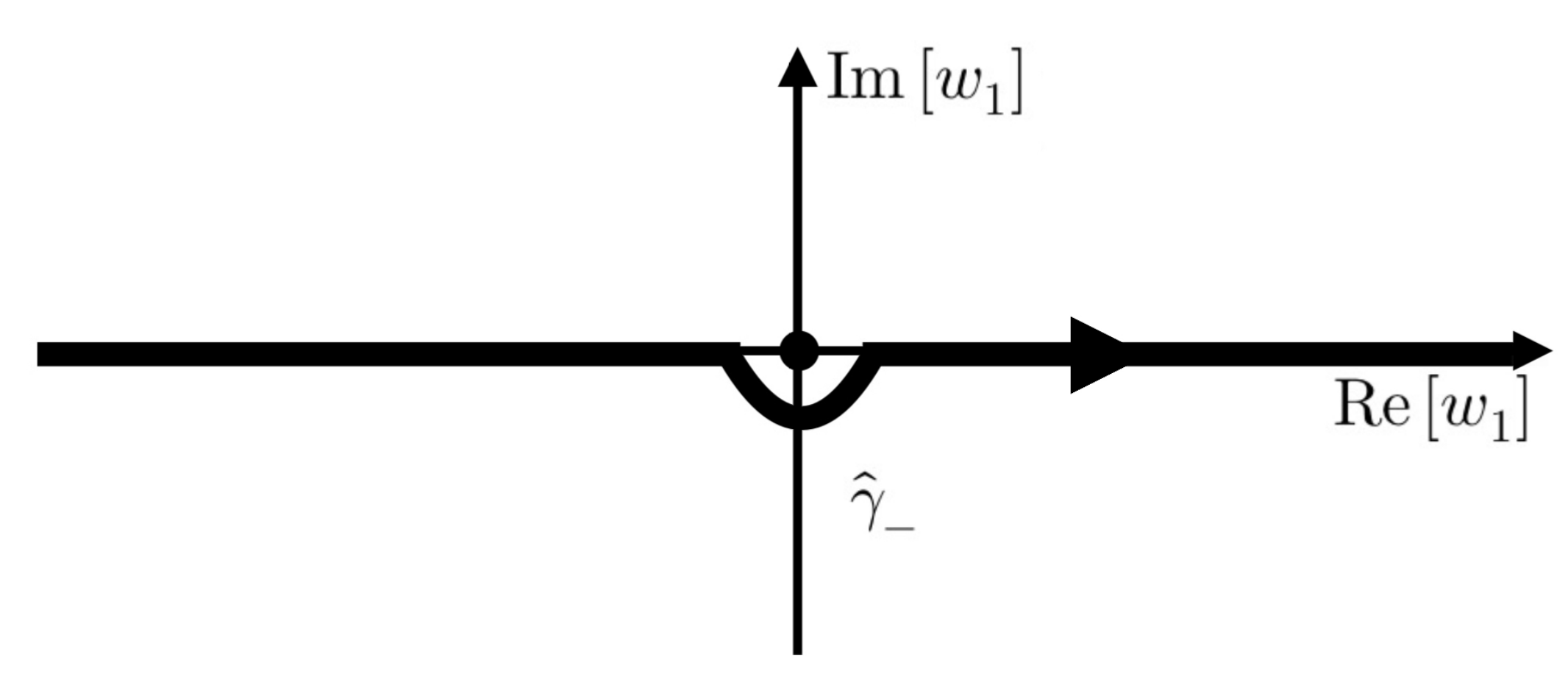}}
\caption{Contours $\hat \gamma_+$ and $\hat \gamma_-$
}
\label{fig01}
\end{figure}

One can construct the fragment of $\Gamma$ in $D$ in one of the two ways: as
\[
\Gamma_+ = \hat \gamma_+ \times [- \varepsilon, \varepsilon] \times [-\varepsilon , \varepsilon]
\]
or as
\[
\Gamma_- = \hat \gamma_- \times [- \varepsilon, \varepsilon] \times [-\varepsilon , \varepsilon].
\]
Here $\hat \gamma_\pm$ are contours that bypass the point $w_1 = 0$ from above and below (see \ \figurename~\ref{fig01}), and $[- \varepsilon, \varepsilon]$ are small fragments of the real axes of the variables $w_2$ and~$w_3$. In the products above, it is assumed that the first contour refers to the variable $w_1$, the second to the variable $w_2$, and the third to the variable~$w_3$.

It can be shown (this is a nontrivial statement, which we leave here without a proof) that any other domain not intersecting $\sigma_j$ and coinciding with $\mathbb{R}^3$ outside the neighborhood of $\sigma_j$ can be locally deformed into $\Gamma_+$ or~$\Gamma_-$.    
This means that there are only two ways to bypass the singularity $\sigma_j$ by the integration domain.

These surfaces stand for the  bypasses ``above'' and ``below'' from one-dimensional complex analysis. It is correct to say that $\Gamma_+$ bypasses above  the singularity $\sigma_j$ {\em in variables $(w_1, w_2, w_3)$}, and $\Gamma_-$ bypasses below it {\em in these variables\/}.  It is assumed that the singularity
corresponds to $w_n = 0$ for one of the variables. 

If $\Gamma$ is described
by a real small vector function $\eta(\xi) = (\eta_1 , \eta_2, \eta_3)$ (see \ (\ref{e:009})), 
the following statement is fulfilled:

\begin{proposition}
Let $D$ be a small ball in $\mathbb{R}^3$ centered at $\xi^* \in \sigma_j$, and let $\Gamma$ be described by a small enough smooth 
vector field~$\eta$.
The following condition is sufficient for $\Gamma \cap D$ not to intersect $\sigma_j$:
\begin{equation}
\eta(\xi) \cdot \nabla g(\xi) \ne 0 , 
\qquad 
\xi \in D.
\label{e:010}
\end{equation}
\end{proposition}

The proof of this proposition is that for small enough $\eta$ one can use the first approximation of the Taylor series
\[
{\rm Im}[g_j (\xi + i \eta (\xi))] \approx
\eta (\xi) \cdot \nabla g_j(\xi)  \ne 0,
\]
which means, 
\[
\xi + i \eta (\xi) \notin \sigma_j.
\]

The condition (\ref{e:010}) means that the vector $\eta(\xi)|_{\sigma'_j}$ 
is nowhere tangential to the surface~$\sigma'_j$.
Hence, the vector
$\eta(\xi)$, $\xi \in \sigma_j'$ points either to one side or to the other side of the surface
$\sigma'_j$, and this side is the same at all points of $\sigma'_j$. 
The made above classification into $\Gamma_+$ and $\Gamma_-$ corresponds to whether the vector $\eta$ is looking upward or downward from
$\sigma'_j$ along the variable~$w_1$.

Thus, to determine the integration domain in (\ref{e:001})  it is sufficient to specify, for
each component of the singularity $\sigma_j$, the side in which the vector $\eta$ looks with respect to the surface~$\sigma'_j$. 
The specification of this side is, of course, not arbitrary, but reflects the physical peculiarities of the problem. We see two ways of defining this direction:

\begin{itemize}

\item
The integral (\ref{e:001}) with the phase function (\ref{e:006}) can arise as the inverse of the Fourier--Laplace transform (for some of the variables $z_m$ the Fourier transform is applied, and for some of the variables the Laplace transform is applied). The inverse of the Laplace transform is given by the Mellin transform. The Mellin transform contains integration along a straight line that is a shifted real axis of~$\xi_m$. 

\item
In the study of wave phenomena, the {\em principle of limiting absorption} is a standard technique. The wave number $k$ in the medium is assumed to have a small imaginary part that corresponds to a small absorption in the medium. An exponentially decaying solution is sought and then a limit is taken at which the imaginary part of $k$ tends to zero. In the language of integrals like (\ref{e:001}), this corresponds to the fact that the functions $g_j$ depend on ${\rm Im}[k]$ as a parameter:
\[
g_j (\xi) = g_j (\xi ; {\rm Im}[k]) 
\]
These functions have the real property only for the limit ${\rm Im}[k] \to 0$.
Considering the case of small nonzero ${\rm Im}[k]$ allows us to determine how $\Gamma = \mathbb{R}^3$ bypasses $\sigma_j$. This side is also preserved in the limit ${\rm Im}[k] = 0$ when the volume $\Gamma$ has to be displaced from~$\mathbb{R}^3$.   
\end{itemize}


\section{The principle of locality}

The principle of locality in our case means the following: 

{\em 
For almost all values of $z$, 
the asymptotic evaluation of the integral (\ref{e:001}) at $\Lambda \to \infty$ is	given	by	the	sum	of	several	terms	defined	by	the
neighborhoods of several special points:
}
\begin{equation}
u(\Lambda) = \sum_\nu u_\nu (\Lambda) + O(e^{- \zeta \Lambda}) 
\label{e:011}
\end{equation}
{\em for some $\zeta > 0$.}

Each of the terms $u_\nu$ is expanded into a series of inverse powers of~$\Lambda$. We will be interested only in the leading term of each of these series.

For us, the principle of locality follows from the statement below, which we will formulate but not prove (the proof is constructive but rather lengthy):

\begin{proposition}
If everywhere, except in the neighborhood of several special points, one can construct the desired and allowed deformation of the domain $\Gamma$ locally (i.e., in small neighborhoods), then everywhere outside the neighborhoods of the special points the desired and allowed deformation can be constructed globally. In other words, local deformations are "glued" into a global deformation.
\end{proposition}

This theorem allows us to restrict ourselves to local deformations $\Gamma$ in~$\Gamma'$. By local deformations we mean  {\em locally constant\/} deformations, i.e.\ such $\eta(\xi; \tau)$, $\tau \in [0,1]$, that the vector $\eta$ is constant along the coordinates $\xi$ in a small neighborhood of the chosen point.


\section{Points that do not contribute to the integral evaluation
}

Here we describe the points that are not special. i.e.\ where a desired and allowed deformation is possible.
Everywhere we will assume that ${\rm Im}[G(\xi^*)] = 0$ (otherwise no deformation of $\Gamma$ is required).
Moreover, we will assume that ${\rm Im}[G(\xi)] = 0$ in some real neighborhood of  $\xi^*$, i.e.\ that $\nabla G$ in this neighborhood is a real vector.

\vskip 6pt
\noindent
{\bf Regular point. }
First of all, consider a point  $\xi^*$, which does not belong to singularities $\sigma_j$,  and in which
\[
\nabla G \ne {\bf 0}.
\]
In this case, the criterion of desirability of deformation is the inequality: 
\begin{equation}
\eta(\xi^* ; 1) \cdot \nabla G(\xi^* ) > 0.
\label{e:012}
\end{equation}
Obviously, for a small neighborhood
$\xi^*$ there exists a desired and allowed deformation for which
\[
\eta(\xi^* ; 1) = \varepsilon \nabla G(\xi^* ), \qquad \varepsilon >0 
\]
for the whole neighborhood of $\xi^*$ ($\varepsilon$ is small). The deformation is given, for example, by a linear interpolation between the initial state 
$\eta(\xi^* ; 0)$ and~$\eta(\xi^* ;1)$.  Thus, in the neighborhood of a regular point, there exists an allowed and desired deformation, and hence the neighborhoods of regular points do not contribute to the evaluation of the integral.

\vskip 6pt
\noindent
{\bf Regular point on the singularity. }
Let $\xi^*$  belong to one of the components of the singularity $\sigma_j$. 
Let 
$\nabla G \ne {\bf 0}$ and 
$ \nabla G \ne c \nabla g_j $
for any scalar value of $c$,  i.e., the vectors $\nabla G$ and $\nabla g_j$ are noncollinear. 
The latter condition can be formulated as follows:
 
\begin{equation}
\nabla \left( 
\left.
G (\xi^* )
\right|_{\sigma'_j}
\right) \ne {\bf 0}.
\label{e:013a}
\end{equation}
(the gradient of the restriction of $G$ on the surface $\sigma'_j$ is non-zero). It turns out that 
in this case it is also always possible to construct the allowed and desired deformation.

Consider in the real space the real vectors
 $\nabla G(\xi^* )$, $\nabla g_j (\xi^*)$ and 
$\eta(\xi^* ; 0)$. $\Gamma$ does not intersect $\sigma_j$, and we assume that this is provided with the condition 
\[
\eta(\xi^*; 0) \cdot \nabla g_j (\xi^*) \ne 0
\]
(see\ (\ref{e:010})).

If the inequality 
$\eta(\xi^* ; 0) \cdot \nabla G(\xi^* ) > 0,$ is fulfilled, then in the neighborhood of $\xi^*$ 
the initial $\Gamma$ is desired, and no 
deformation is required. 

Let be $\eta(\xi^* ; 0) \cdot \nabla G(\xi^* ) < 0$. 
Construct the vector
\[
{\bf a} = \nabla  g_j (\xi^*) \times (  \nabla  g_j (\xi^*) \times \nabla  G(\xi^* ))
\]
in $\mathbb{R}^3$. Note that
\begin{equation}
{\bf a} \cdot \nabla  g_j (\xi^*) = {\bf 0},
\qquad
{\bf a} \cdot \nabla  G (\xi^* ) \ne {\bf 0}
\label{e:013}
\end{equation}
(the second part requires application of the well-known identity about vector and scalar product).

The sought deformation of the domain $\Gamma$ in the neighborhood of $\xi^*$ is given by the formula
\begin{equation}
\eta({\rm Re}[\xi] ; \tau) = \eta({\rm Re}[\xi] ; 0)  + \alpha \tau {\bf a}.
\label{e:014a}
\end{equation}
at a suitable choice of sign and value of $\alpha$. The first part of (\ref{e:013}) says that moving along ${\bf a}$ in the first approximation does not lead to an intersection of $\sigma_j$, and the second part says that such a movement can lead to a change of sign of $\eta \cdot \nabla G(\xi^* )$. For example, we can choose
\[
\alpha = - 2 \frac{\eta(\xi^* ; 0) \cdot \nabla G(\xi^* )}{{\bf a} \cdot \nabla  G (\xi^* )} .
\]

Thus, we come to an important conclusion: a point belonging to a single singularity $\sigma_j$ under the condition
(\ref{e:013a}) does not yield a contribution to the asymptotics. This differs from the one-dimensional case, where singular points usually give terms in asymptotics (see also \cite{Mironov2021,Assier2022}).

\vskip 6pt
\noindent
{\bf The regular point on a crossing point of two singularities.}
Let $\xi^* \in  \sigma_j \cap \sigma_m$ and the inequality is fulfilled 
\begin{equation}
\nabla \left( 
\left.
G(\xi^* )
\right|_{\sigma'_j \cap \sigma'_m}
\right) \ne {\bf 0}.
\label{e:014}
\end{equation}
Note that the crossing of $\sigma'_j \cap \sigma'_m$ is a one-dimensional set
(line) in three-dimensional domain, so that the gradient in (\ref{e:014}) is a derivative along this line. Assume that the crossing of singularities is transversal (the corresponding surfaces do not touch).

In this case, the allowed and desired deformation is constructed similarly to  (\ref{e:014a}), and the vector ${\bf a}$ is defined as
\[
{\bf a} = \nabla g_j (\xi^*) \times \nabla g_m (\xi^*). 
\]

This result distinguishes the three-dimensional case from the two-dimensional case, since in the two-dimensional case the intersections of singularities, as a rule, lead to the appearance of terms in the asymptotics.


\section{Points contributing to the integral asymptotics (special points)
}

Here we list the main types of special points. The classification is not exhaustive, since, for example, there can exist  degenerate cases (some derivatives or determinants may vanish). Typically, this corresponds to a merge of special points of the listed types. We assume, however, that consideration of the degenerate cases causes no problem.

\begin{figure}[ht]
\centerline{\epsfig{width = 10 cm, file=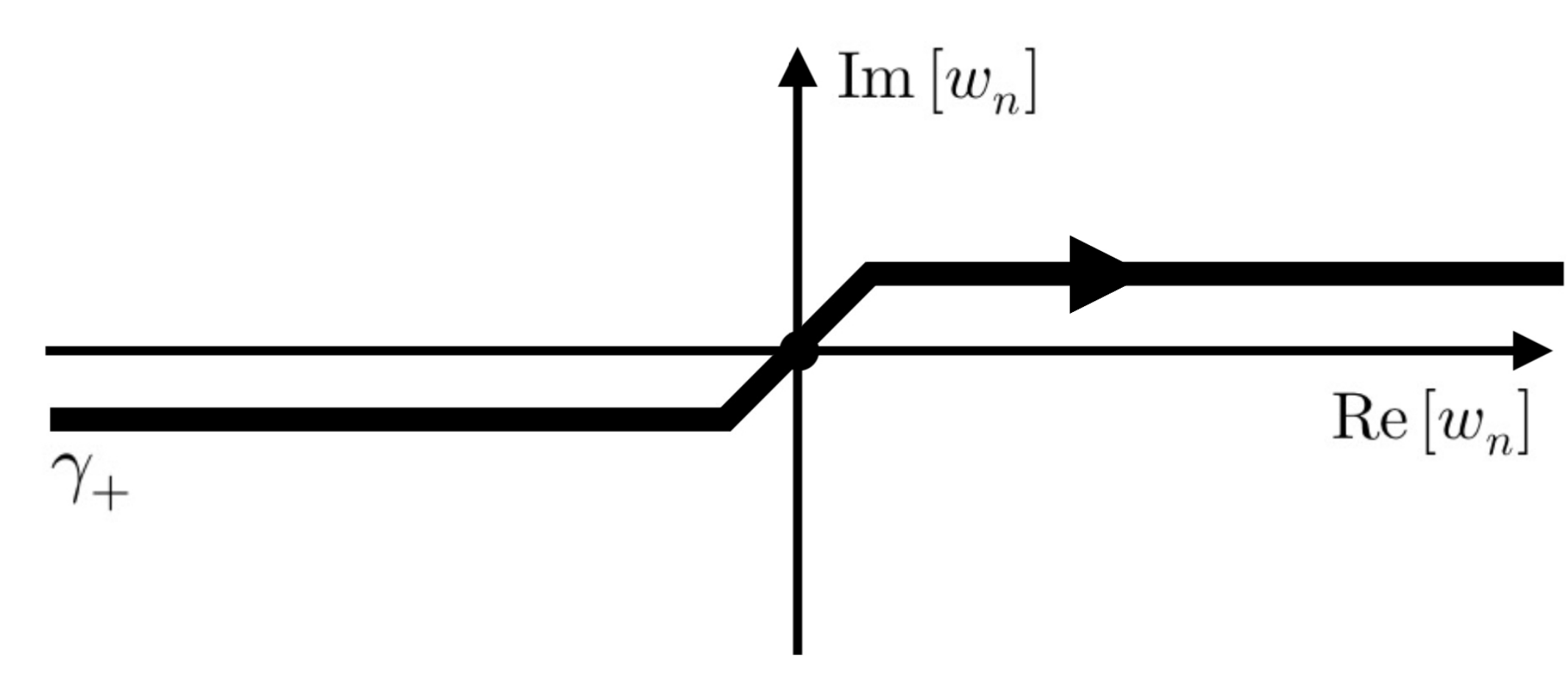}}
\centerline{\epsfig{width = 10 cm, file=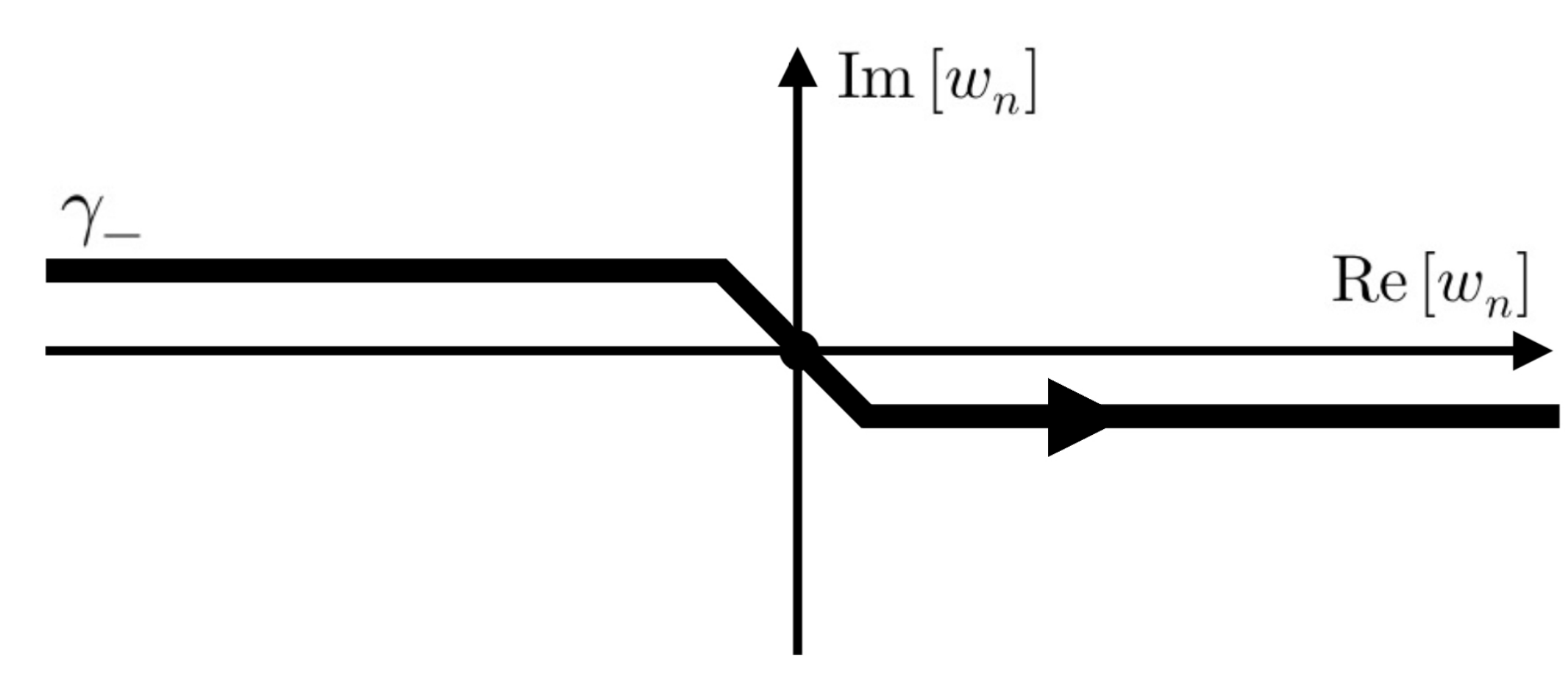}}
\caption{Contours $\gamma_+$ and $\gamma_-$ for the stationary phase point
}
\label{fig02}
\end{figure}

\vskip 6pt
\noindent
{\bf Stationary phase point.}
Let at the point $\xi^* \in \mathbb{R}^3$ the functions $G$ and $F$ be holomorphic,  
\[
\nabla G(\xi^* ) = {\bf 0},
\]
and the determinant of the Hessian
\[
H \equiv \left( \begin{array}{ccc}
\frac{\ptl^2 G}{\ptl \xi_1^2} & \frac{\ptl^2 G}{\ptl \xi_1 \ptl \xi_2} & \frac{\ptl^2 G}{\ptl \xi_1 \ptl \xi_3} \\
\frac{\ptl^2 G}{\ptl \xi_1 \ptl \xi_2} & \frac{\ptl^2 G}{\ptl  \xi_2^2} & \frac{\ptl^2 G}{\ptl \xi_2 \ptl \xi_3} \\
\frac{\ptl^2 G}{\ptl \xi_1 \ptl \xi_3} & \frac{\ptl^2 G}{\ptl \xi_2 \ptl \xi_3} & \frac{\ptl^2 G}{\ptl \xi_3^2} 
\end{array} \right) 
\]
at the point $\xi^*$ be non-zero. Such a point is the usual three-dimensional stationary phase point or a saddle point
\cite{Fedoryuk1977}. We write out the consideration here for completeness. Let us bring the Hessian to the principal axes
by going to the local 
variables $(w_1, w_2, w_3)$ obtained from $\xi$ by a real nondegenerate linear transformation. In the new variables
\[
G(w ) = G(\xi^* ) + \frac{1}{2} \sum_{n = 1}^3 \beta_n w_n^2 + O(|w|^3) 
\]
with some real non-zero constants $\beta_n$. 
Then the local deformation of the domain $\Gamma$ is given in new variables by the formula
\[
\Gamma' = \gamma_{{\rm sign} (\beta_1)} \times \gamma_{{\rm sign} (\beta_2)} \times \gamma_{{\rm sign} (\beta_3)},
\]
where the contours $\gamma_\pm$ are shown in \figurename~\ref{fig02}.
Obviously, this deformation is allowed, and it is desired everywhere except at the point~$\xi^*$. 

If the function $F(\xi)$ is not zero at the point $\xi^*$, the leading term of the asymptotics of the corresponding contribution $u_\nu$ can be easily obtained: 
\begin{equation}
u_\nu (\Lambda) = F(\xi^*)\exp \left\{i \Lambda G(\xi^*) + i \frac{\pi}{4}\sum_{n= 1}^3 
{\rm sign}(\beta_n)
\right\}  \frac{(2\pi)^{3/2} J}{\sqrt{|\beta_1 \beta_2 \beta_3|} \Lambda^{3/2}} + O(\Lambda^{-5/2}),
\label{e:014b}
\end{equation}
$J$ (here and below) is the determinant of the Jacobian of the transformation
\begin{equation}
J = {\rm det}\left( 
\frac{\ptl (\xi_1 , \xi_2 , \xi_3)}{\ptl (w_1 , w_2 , w_3)}
\right) .
\label{e:014c}
\end{equation}

\vskip 6pt
\noindent
{\bf Stationary phase point on a single singularity.}
Let $\xi^*$ belong to a single component of the singularity $\sigma'_j$, it is the stationary phase point of the restriction of 
$G$ to $\sigma'_j$:
\begin{equation}
\nabla \left( 
\left.
G (\xi^*)
\right|_{\sigma'_j}
\right) = {\bf 0},
\label{e:015}
\end{equation}
and the Hessian matrix of the restriction of $G$ to  $\sigma'_j$ is nondegenerate. Assume also
that $\xi^*$ is not a stationary phase point, i.e.\ $\nabla G(\xi^*) \ne {\bf 0}$. 

Perform near $\xi^*$ a biholomorphic real coordinate change $(\xi_1, \xi_2, \xi_3) \to (w_1, w_2, w_3)$ such that 
\begin{equation}
w_1 = \alpha g_j(z),
\label{e:015z}
\end{equation}
\begin{equation}
G(w) = G(\xi^* ) + w_1 + \frac{1}{2}( \beta_2 w_2^2 + \beta_3 w_3^2 ) + \dots
\label{e:015n}
\end{equation}
where $\alpha$, $\beta_1$, $\beta_2$ are non-zero real coefficients. By biholomorphic we mean
such a coordinate change that the functions 
$w_j (\xi)$ and $\xi_j (w)$, $(w \equiv (w_1, w_2, w_3))$, are holomorphic in the domain under consideration. We assume that the substitution is  {\em real}, i.e.\ real $\xi$ correspond to real $w$ and vice versa. By choosing $\alpha$ we ensure that the coefficient at $w_1$ in (\ref{e:015n}) is equal to unity.

\begin{figure}[ht]
\centerline{\epsfig{width = 10 cm, file=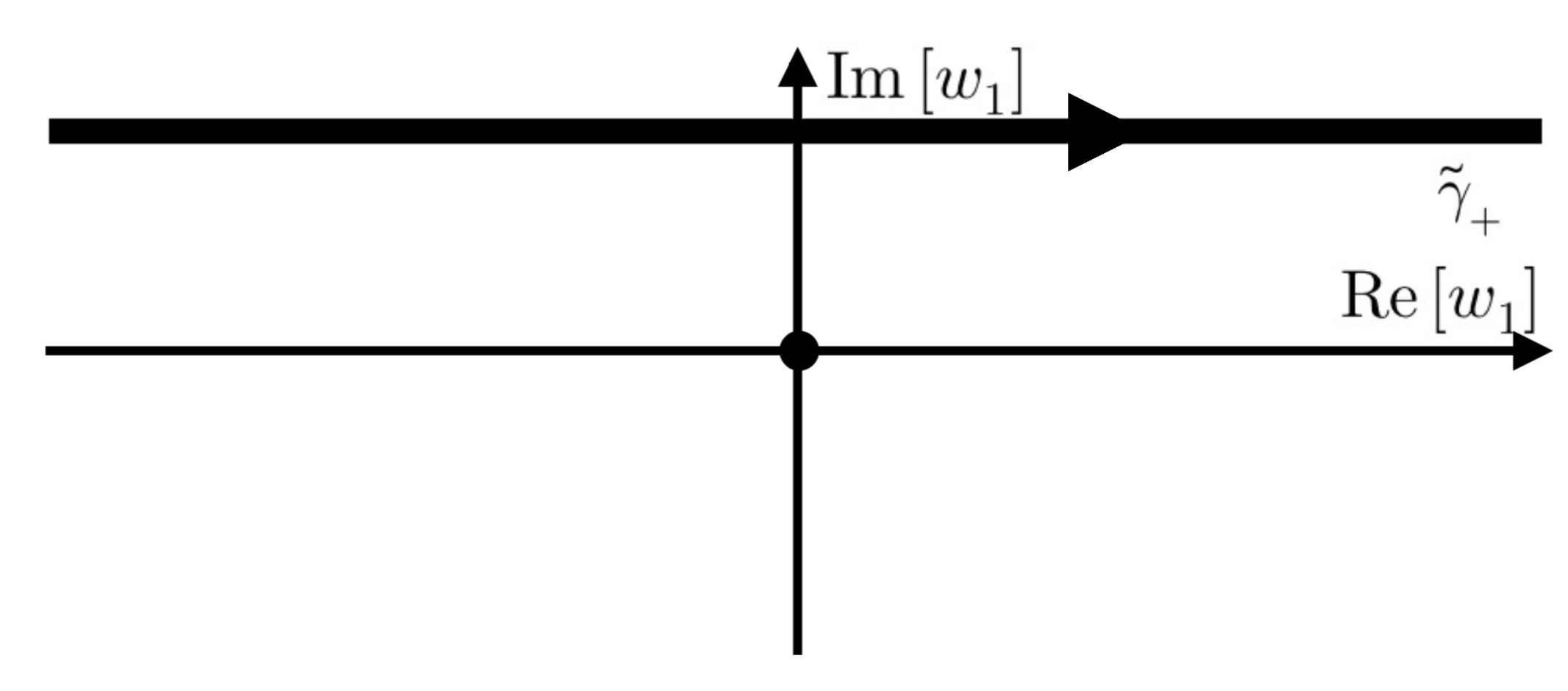}}
\centerline{\epsfig{width = 10 cm, file=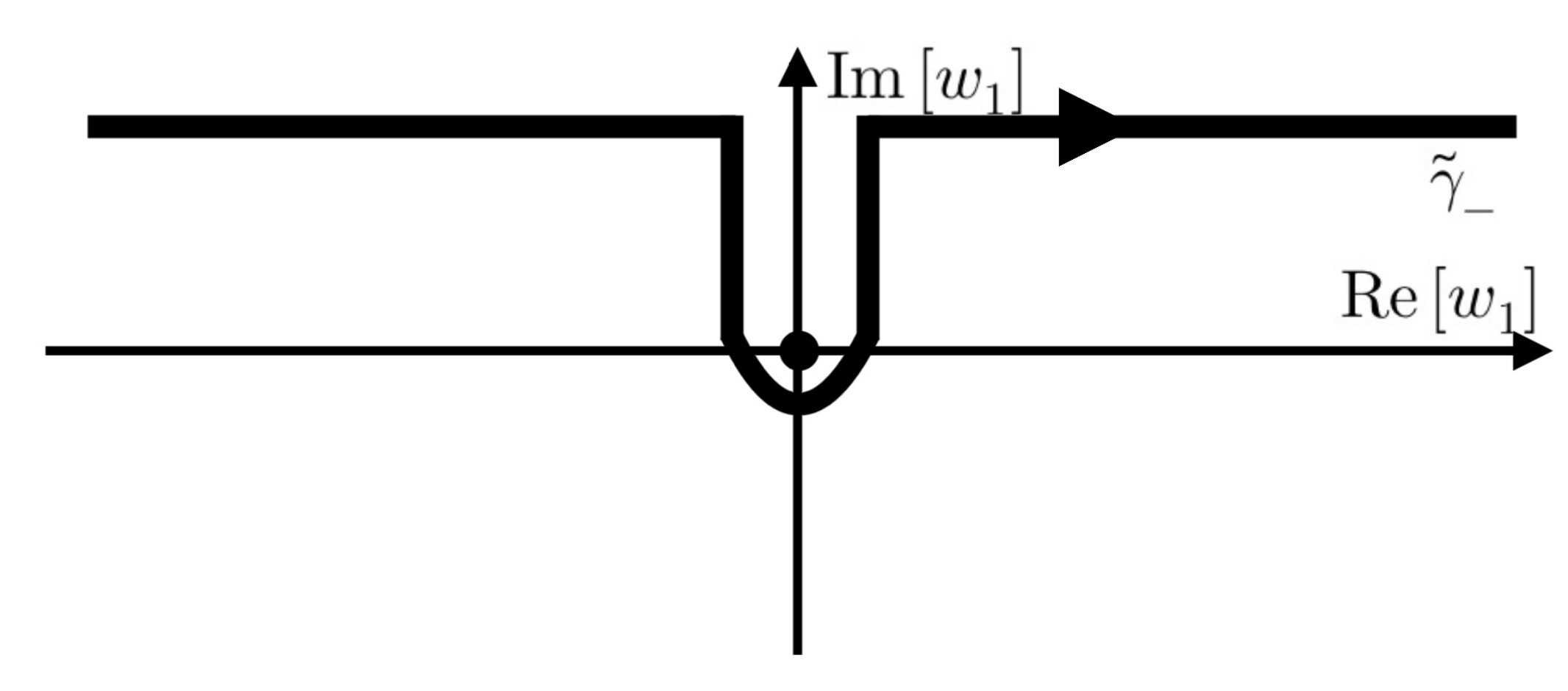}}
\caption{Contours $\tilde \gamma_+$ and $\tilde \gamma_-$ for the stationary phase point
}
\label{fig03}
\end{figure}

The singularity is given by $w_1 = 0$ in the new coordinates. The original integration domain passes 
above or below the set 
$w_1 = 0$. 
The deformed domain $\Gamma'$ is constructed near $\xi^*$ 
in coordinates $(w_1 , w_2 , w_3)$ as 
\[
\Gamma' = \tilde \gamma_{\pm} \times \gamma_{{\rm sign}(\beta_2)} \times \gamma_{{\rm sign}(\beta_3)}.
\]
The contours $\tilde \gamma_{\pm}$ are shown in \figurename~\ref{fig03}. The contour $\tilde \gamma_+$ is chosen if
$\Gamma$ passes above $w_1 = 0$, and $\tilde \gamma_-$ if it passes below.  
In both cases, this deformation is allowed, but in the second case it is not desired in the neighborhood of the point $\xi^*$.
In the first case, the neighborhood of the point $\xi^*$ does not contribute to the asymptotics of the integral, but in the second case it does. The integral in the second case
case can be estimated as a nested one in the coordinates $(w_1, w_2, w_3)$:
\begin{equation}
u_\nu (\Lambda) =  \Phi(\Lambda)\exp \left\{i \Lambda G(\xi^* ) + i \frac{\pi}{4}\sum_{n= 2}^3 
{\rm sign}(\beta_n)
\right\}  \frac{2 \pi J}{\Lambda \sqrt{|\beta_2 \beta_3|}}  ,
\label{e:015a}
\end{equation}
$J$ is the determinant of the Jacobian of the transformation (see\ (\ref{e:014c})), $\Phi$ is the first term of the asymptotic evaluation of the integral 
\begin{equation}
\Phi' (\Lambda) = \int_{\tilde \gamma_-} F (w_1 , 0 ,0) e^{i \Lambda w_1} dw_1
\label{e:015b}
\end{equation}
at large $\Lambda$.
This integral is easily evaluated by standard methods in the case of power singularities of~$F$.
Namely, let the leading term of the three-dimensional Taylor series expansion of $F$ near $\xi^*$ be
\begin{equation}
F(w_1 , w_2, w_3) = C(w_2,w_3) w_1^\mu,
\label{e:015y}
\end{equation}
where $C(w_2,w_3)$ is a function holomorphic at $w_2=0$, $w_3=0$ and not equal to zero there. Here we take only the simple pole or the branch point for the variable $w_1$ (so either $\mu=-1$ or $\mu\notin\mathbb{Z}$) Then
\begin{equation}
\Phi(\Lambda) = C(0,0) \Lambda^{-\mu - 1} I(\mu),
\label{e:015c}
\end{equation}
where 
\[
I(-1)  = 2\pi i,
\]
and
\[
I(\mu) = e^{\pi i (\mu + 1)/2}(1 - e^{-i 2\pi \mu }) \Gamma_E(1 + \mu), \qquad \mu \notin \mathbb{Z}.
\]
($\Gamma_E$ is the Euler Gamma function).

We conclude that the stationary phase point on the singularity contributes to the asymptotics of the integral if some additional condition is fulfilled. This condition can be formulated as follows: near this point, the initial integration domain $\Gamma$ bypasses the singularity $\sigma_j$ from the growth side of the exponential multiplier~$e^{i \Lambda G}$.

\vskip 6pt
\noindent
{\bf Stationary phase point on the line of crossing of two singularities.
}
Let $\xi^*$ belong to the crossing of two singularities: $\xi^* \in \sigma_j' \cap \sigma_m'$,
and
\begin{equation}
\nabla \left( 
\left.
G(\xi^* )
\right|_{\sigma'_j \cap \sigma'_m}
\right) = {\bf 0}.
\label{e:016}
\end{equation}
Let us make a local biholomorphic real variable change $(\xi_1 , \xi_2, \xi_3) \to (w_1 , w_2 , w_3)$ such that
\begin{equation}
w_1 = \alpha_1 g_j(\xi), 
\qquad 
w_2 = \alpha_2 g_m(\xi), 
\label{e:016a}
\end{equation}
for some $\alpha_1$, $\alpha_2$, 
to make $G$ having the form 
\begin{equation}
G = G(\xi^* ) + w_1 + w_2 + \frac{1}{2} \beta w_3^2  + \dots
\label{e:016b}
\end{equation}
for some $\beta$ (we assume that such a substitution is possible, i.e., that $\xi$ is not a stationary phase point for either of the two singularities, that the crossing of the singularities is transversal, and that the second derivative of $G$ along the line of intersection of the singularities is non-zero).

The singularities in the new coordinates are given by $w_1 = 0$ and $w_2 = 0$. The initial domain $\Gamma$ can bypass 
above or below
each of these sets. Thus, $\Gamma$
can be locally deformed into 
\[
\Gamma' = \tilde \gamma_{\pm} \times \tilde \gamma_\pm \times \gamma_{{\rm sign}(\beta)}.
\]
This deformation is desired if at least one of the first two contours is $\tilde \gamma_+$, i.e., if at least one singularity is bypassed from above (in $w$ coordinates).
Thus, the condition that the point $\xi^*$ contributes to the asymptotics of the integral is that $\Gamma$ bypasses below both singularities. The estimation of this contribution in this case is 
\begin{equation}
u_{\nu} (\Lambda) = \Phi(\Lambda)
\exp \left\{i \Lambda G(\xi^* ) + i \frac{\pi}{4} 
{\rm sign}(\beta)
\right\}  \frac{\sqrt{2 \pi} J}{\sqrt{ \Lambda  |\beta|}}  ,
\label{e:017}
\end{equation}
where $J$ is the Jacobian determinant and $\Phi$ is the leading term of the asymptotics of the integral
\begin{equation}
\Phi'(\Lambda) = \int_{\tilde \gamma_-} \int_{\tilde \gamma_-} F(w_1, w_2 , 0) e^{i \Lambda (w_1 + w_2)} dw_1 \, dw_2   
\label{e:018}
\end{equation}
The latter multiplier can be easily constructed, e.g., if leading term of $F$ near $\xi^*$ has the form
\begin{equation}
F(w_1 , w_2 , w_3) \approx C(w_3) w_1^{\mu_1} w_2^{\mu_2}, 
\label{e:019}
\end{equation}
where $C(w_3)$ is a function holomorphic at zero and not equal to zero there. In this case
\begin{equation}
\Phi(\Lambda) = C(0) \Lambda^{-\mu_1-\mu_2 -2} I(\mu_1) I(\mu_2).
\label{e:019a}
\end{equation}

\vskip 6pt
\noindent
{\bf The crossing of three singularities.}
Let $\xi^*$ belong to an
intersection of three singularities: $\xi^* \in \sigma'_j \cap \sigma'_l \cap \sigma'_m$ (note that such 
an intersection is a discrete set of points). Suppose that $\xi^*$ is not a stationary phase point on any of the lines of pairwise crossing of singularities. Then, by a local biholomorphic real substitution
of the form
\[
w_1 = \alpha_n g_j(\xi) , 
\qquad 
w_2 = \alpha_2 g_l(\xi) , 
\qquad 
w_3 = \alpha_3 g_m(\xi) 
\]
we can reduce the function $G$ to 
\[
G(\xi ) = G(\xi^*) + w_1 + w_2 + w_3 + O(|w|^2).
\]

The
singularities $\sigma_j$, $\sigma_l$, $\sigma_m$ have the form $w_1 = 0$, $w_2 = 0$, $w_3 = 0$, respectively. The initial integration domain $\Gamma$ bypasses above or below each of these singularities. In the coordinates $w$, 
$\Gamma$ can be locally deformed by an allowed deformation into
\[
\Gamma' = \tilde \gamma_{\pm} \times \tilde \gamma_{\pm} \times \tilde \gamma_{\pm},
\]
where the sign $\pm$ in the $n$-th position is chosen depending on whether   
$\Gamma$ bypasses above or below $w_n  = 0$.
The deformation is desired if at least one of these contours is $\tilde \gamma_+$. Hence,
the  condition that  $\xi^*$ 
gives a contribution to the asymptotics is as follows: 
in the $w$ coordinates, the volume $\Gamma$ must bypass below all  three singularities. 

Let the leading term $F$ near $\xi^*$ have the form
\[
F(w_1, w_2, w_3) \approx C w_1^{\mu_1} w_2^{\mu_2} w_3^{\mu_3},
\]
where $C$ is a nonzero constant. Then the evaluation of the contribution of the point $\xi^*$ is given by the evaluation of the nested 
integral
\begin{equation}
u_\nu (\Lambda) \approx C e^{i \Lambda G(\xi^*)}\prod_{n = 1}^3 \int_{\tilde \gamma_-} w_n^{\mu_n} 
e^{i \Lambda w_n } dw_n \approx  
C e^{i \Lambda G(\xi^*)}\prod_{n = 1}^3 \Lambda^{-\mu_n -1} I(\mu_n)
.   
\label{e:020}
\end{equation}

\vskip 6pt
\noindent
{\bf Conical point}.
The special points discussed above were generalizations of special points in the one-dimensional and two-dimensional cases. In the 
three-dimensional case, however, fundamentally new conical points can arise.

Let there exist coordinates $(w_1, w_2, w_3)$ in which the function $g_j$ defining the singularity $\sigma_j$ has the form
\begin{equation}
g_j (w_1, w_2, w_3)  = w_1^2 + w_2^2 - w_3^2.
\label{e:021}
\end{equation}
As before, we assume that the substitution $(\xi_1, \xi_2 , \xi_3) \to (w_1, w_2, w_3)$ is biholomorphic and real. The values $w_1 = w_2 = w_3 = 0$ correspond to the point $\xi = \xi^*$. Obviously, $\sigma'_j$ represents a double-sided cone.  

We also assume that $G$ is approximated by a linear function
\[
G(\xi) = G(\xi^*) + \alpha_1 w_1 + \alpha_2 w_2 + \alpha_3 w_3 + O(|w|^2)
\]
with a non-zero vector $(\alpha_1, \alpha_2, \alpha_3)$.

A conical singular point usually occurs when the problem considers wave propagation in two spatial and one temporal coordinates, and at low spatial and temporal frequencies the wave velocity is non-zero and finite. In this case, $g_j = 0$ is the dispersion equation for the medium.

First, let us describe how the initial integration domain $\Gamma$ can be located near the conical singular point and not intersect~$\sigma_j$. Let us assume that in the coordinates $w$ near the origin
\[
\Gamma = \mathbb{R}^3 + i (\epsilon_1, \epsilon_2, \epsilon_3),
\]
i.e.\ $\Gamma$ 
is $\mathbb{R}^3$ shifted by a small constant imaginary vector. 
Let some point of $\Gamma$ be $w = (w_1' + i \epsilon_1, w_2' + i \epsilon_2 , w_3' + i \epsilon_3)$, 
where $w'_{1,2,3}$ are real.
Let us write out the real and imaginary parts of $g_j$ and write down the conditions that $w \in \sigma_j$:
\begin{equation}
(w_1')^2 + (w_2')^2 - (w_3')^2 = \epsilon_1^2 + \epsilon_2^2 - \epsilon_3^2,
\label{e:022}
\end{equation}
\begin{equation}
w_1' \epsilon_1 + w'_2 \epsilon_2 - w'_3 \epsilon_3 = 0.
\label{e:023}
\end{equation}
This system of equations has no real solutions  $w'$ if and only if 
\begin{equation}
\epsilon_3^2 > \epsilon_1^2 + \epsilon_2^2.
\label{e:024}
\end{equation}
One can see that  (\ref{e:022}) is a two-sheet hyperboloid lying inside the cone $(w_1')^2 + (w_2')^2 - (w_3')^2 = 0$, and the plane (\ref{e:023}) intersects this cone
only at the origin. Thus, (\ref{e:024}) describes possible integration domains near the singular point.

\begin{figure}[ht]
\centerline{\epsfig{width = 10 cm, file=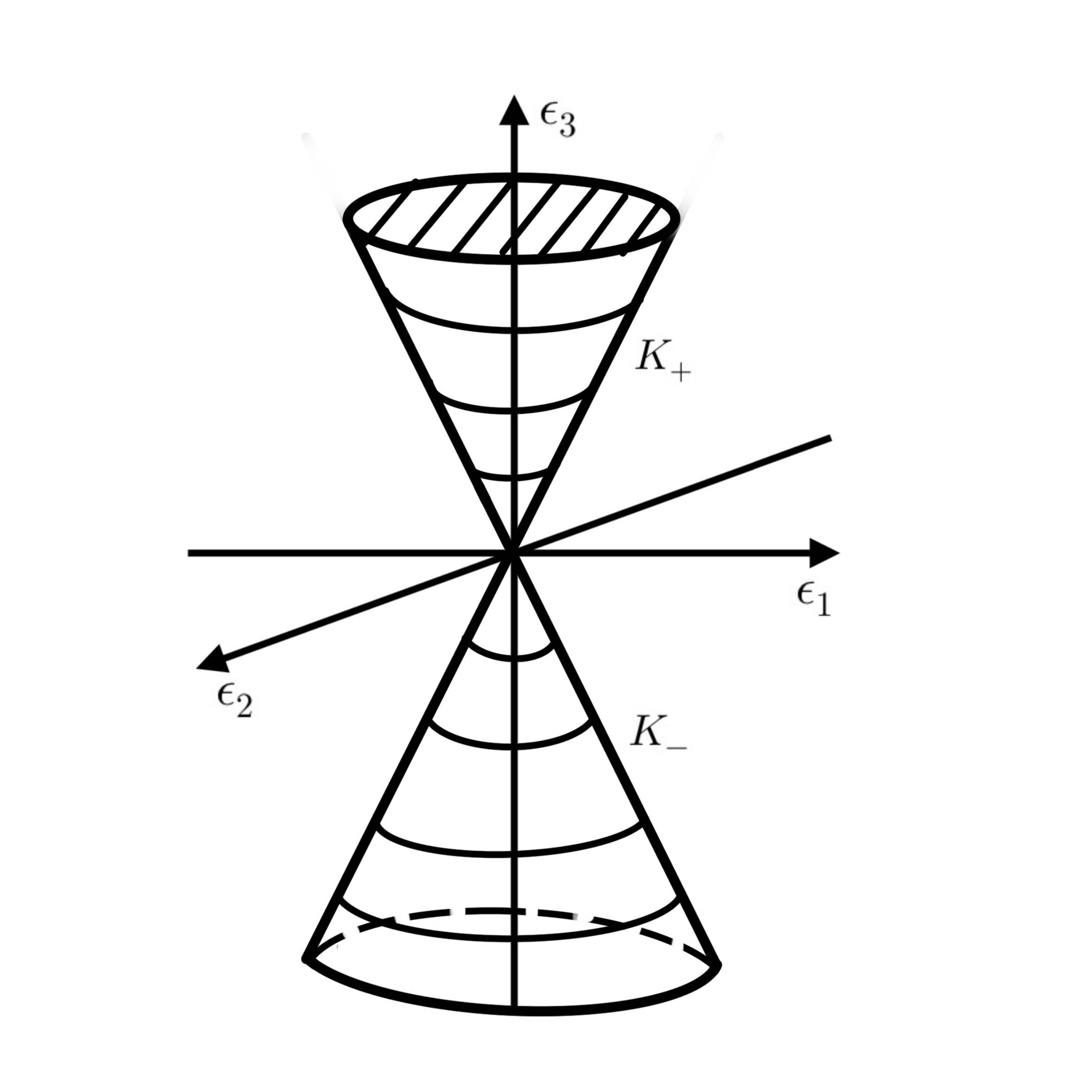}}
\caption{Cones $K_\pm$}
\label{fig07}
\end{figure}

The set (\ref{e:024}) is the union of two one-sided cones
 $K_\pm$:
\[
K_+ : \quad \epsilon_3 > \sqrt{\epsilon_1^2 + \epsilon_2^2},
\]
\[
K_-: \quad\epsilon_3<-\sqrt{\epsilon_1^2 + \epsilon_2^2}.
\]
These cones are shown in \figurename~\ref{fig07}.
The allowed deformations of the domain $\Gamma$ near the conical point consist of moving the point $(\epsilon_1, \epsilon_2, \epsilon_3)$ within one of the two cones~$K_\pm$.
One cannot move from one cone to the other
without crossing $\sigma_j$.
Thus, there are two ways to bypass a singularity $\sigma_j$ with a conical point: 
informally speaking, $(\epsilon_1, \epsilon_2 , \epsilon_3) \in K_+$ corresponds to
bypassing above the singularity, and 
$(\epsilon_1, \epsilon_2 , \epsilon_3) \in K_-$ corresponds to bypassing below the singularity. 

Let the deformation $\Gamma \to \Gamma'$ correspond to a change of the displacement
vector from
$(\epsilon_1, \epsilon_2 , \epsilon_3)$ to $(\epsilon'_1, \epsilon'_2 , \epsilon'_3)$
The deformation is desired if 
\begin{equation}
\epsilon_1' \alpha_1 + \epsilon_2' \alpha_2 + \epsilon_3' \alpha_3 > 0.
\label{e:025}
\end{equation}
Obviously, for
\[
|\alpha_3| < \sqrt{\alpha_1^2 + \alpha_2^2}
\]
points satisfying (\ref{e:025}) exist both in $K_+$, and in $K_-$, i.e.\ the desired and allowed deformation exists for any 
mutual position of $\Gamma$ and~$\sigma_j$. 
The condition that the desired and allowed deformation does not exist is the following: 
\begin{equation}
\epsilon_3 < - \sqrt{\epsilon_1^2 + \epsilon_2^2} 
\qquad 
\mbox{and}
\qquad 
\alpha_3 > \sqrt{\alpha_1^2 + \alpha_2^2} 
\label{e:026}
\end{equation}
or
\begin{equation}
\epsilon_3 > \sqrt{\epsilon_1^2 + \epsilon_2^2} 
\qquad 
\mbox{and}
\qquad 
\alpha_3 < -\sqrt{\alpha_1^2 + \alpha_2^2} .
\label{e:027}
\end{equation}
If either of these two conditions is fulfilled, the point $\xi^*$ gives a contribution to the integral asymptotics.

Let condition (\ref{e:026}) be fulfilled. Let us demonstrate an example of calculation of the asymptotics of the integral 
in the neighborhood of a conical point. 
Without attempting to consider the general case, we restrict ourselves to the situation when the higher term in the power expansion of $F$ near $\xi^*$ is given by the formula 
\begin{equation}
F(w_1, w_2, w_3) = \frac{C}{w_1^2 + w_2^2 - w_3^2} + O(|w|^3).
\label{e:028}
\end{equation}
The leading term of the estimate is 
\begin{equation}
u_{\nu}(\Lambda) = C\exp \{i \Lambda G(\xi^*) \} J 
\int \limits_{-\infty - i\epsilon}^{\infty - i \epsilon}
\int \limits_{-\infty }^{\infty }
\int \limits_{-\infty }^{\infty }
\frac{e^{i \Lambda( \alpha_1 w_1 + \alpha_2 w_2 + \alpha_3 w_3)}}{w_1^2 + w_2^2 - w_3^2}
dw_1 \, dw_2 \, dw_3
\label{e:029}
\end{equation}
The triple integral is easily taken:
\begin{equation}
u_{\nu}(\Lambda) = \frac{4C\pi^2 J \exp \{i \Lambda G(\xi^*) \} }{\Lambda \sqrt{\alpha^2_3 - \alpha_1^2 - \alpha_2^2}}
 .
\label{e:030}
\end{equation}



\section{Special points on the boundary of the integration domain}

We assumed above that the integration domain $\Gamma$ is unbounded and is a slightly deformed real domain $\mathbb{R}^3$.
Sometimes, it is necessary\footnote{The importance of this issue was brought to the attention of the authors by Prof.\ Daan Huybrechs, see\ \cite{Huybrechs2007}.} to consider also bounded domains. For example, the integration domain may be bounded by a sphere, several planes, etc. Asymptotic terms caused by the boundary may arise in this case.

Let the integration domain be bounded by several surfaces $s_j$ belonging to analytic sets 
$S_j$: $s_j \subset S_j \subset \mathbb{C}^3$. 
Each of the analytic sets is described by an equation
\begin{equation}
S_j : \qquad h_j(\xi) = 0,
\label{e:031}
\end{equation}
where $h_j(\xi)$ is a holomorphic function. We assume that the functions $h_j$ also have the real property introduced earlier. Examples of such analytic sets are the complexified sphere
\begin{equation}
S_j : \qquad \xi_1^2 + \xi_2^2 + \xi_3^2 - R^2 = 0
\label{e:031a}
\end{equation}
or the complexified plane
\begin{equation}
S_j : \qquad \alpha_1 \xi_1 +\alpha_2 \xi_2 + \alpha_3 \xi_3  -\alpha_0 = 0.
\label{e:031b}
\end{equation}
The corresponding elements of the boundary $s_j$ can be a real sphere and a real plane.  
The sets $S_j$ have real dimension 4, and the sets $s_j$ have real dimension 2.
The sets $S_j$ are {\em complexifications\/} of the elements of the boundary~$s_j$.

The Cauchy theorem guarantees the invariance of the integral (\ref{e:001}) under deformation of the integration domain if singularities are not affected and the boundary of the volume remains in place. It turns out that in 
multidimensional complex analysis the condition of invariance of the boundary can be relaxed. 
That is, it is possible to deform the boundary during deformation, but in such a way that elements
$s_j$ are deformed within their complexifications~$S_j$. This statement is based on the fact that
\begin{equation}
\int_{\Gamma''}  Q(\xi) \, d\xi = 0,
\label{e:032}
\end{equation}
if the function $Q(\xi)$ is holomorphic in a neighborhood of $\Gamma''$, and the region $\Gamma''$ lies in some analytic set~$S$. In our case, indeed,
\[
Q(\xi) = F(\xi) e^{i \Lambda G(\xi)},
\]
$S$ is one of the sets $S_j$, and $\Gamma''$ is the region spotted by the boundary part $s_j$ when it is deformed. The statement (\ref{e:032}) can be easily illustrated for the case 
\[
S_j: \quad \xi_1 = {\rm const}.
\]
Namely, within $S_j$ we have $d \xi_1 = 0$, and hence the volume element for $\Gamma''$ is zero:
\[
d\xi |_{S_j} = 0.
\]
Any other case can be reduced to (\ref{e:031b}) by a bigolomorphic coordinate change. 

Then we can proceed according to the scheme used above: we introduce
the concept of the desired deformation of the boundary (i.e., the deformation that makes the integrand on the boundary exponentially small), and then the concept of special points as the points of the boundary where the desired deformation is not possible.

By analogy to $\sigma_j'$ we introduce the real trace $S_j$ as 
\[
S_j' = S_j \cap \mathbb{R}^3.
\]
It is not difficult to see that the special points on the boundary are:
\begin{itemize}
\item
stationary phase points on one element of the boundary $S'_j$, 

\item
stationary phase points at the crossing of boundary elements $S'_j \cap S'_m$, 

\item
points of triple crossing of the boundaries     $S'_j \cap S'_m \cap S'_n$,

\item
conical boundary points.
\end{itemize}

Besides, we can consider various crossings of boundaries and singularities and saddle points for them.

We do not develop these ideas in this paper, but the corresponding formulas for the asymptotics can be constructed easily.

\section{An example: Kelvin waves}
\vskip 6pt
\noindent
{\bf Derivation of the integral for Kelvin waves.} Consider a classical example: evaluation of the complex three-dimensional Fourier integral for Kelvin waves. Let an incompressible fluid of density $\rho$
 occupy in equilibrium a half-space
$(x_1, x_2, x_3)$, $x_3 < 0$. A gravitational force with the free-fall acceleration~$a$ acts in the negative direction of~$x_3$. A small body of volume $V$ is located under the surface of the fluid at the point $(0,0,0)$
(for convenience we assume that it is a cylinder of the length $l$ and cross-section $A$ oriented along the $x_1$-axis). At time $t = 0$ this body starts to move in the positive direction of the axis $x_1$ 
with velocity~$v$. It is necessary to find the vertical displacement $x_3 = H(x_1, x_2, t)$ of the fluid boundary.

The solution of the problem is described by the integral 
\begin{equation}
H(x_1, x_2, t) = \frac{i v V}{8 \pi^3} \int \limits_{-\infty + i \epsilon}^{\infty + i \epsilon}
\int \limits_{-\infty}^{\infty }
\int \limits_{-\infty}^{\infty }
\frac{\omega k_1 e^{i k_1 x_1 + i k_2 x_2 - i \omega t}}{(\omega - k_1 v)(\omega^2 - a \sqrt{k_1^2 + k_2^2})}
dk_1 \, dk_2  \, d\omega.
\label{e:101}
\end{equation}

The derivation of the integral representation (\ref{e:101}) is as follows. 
Let us describe the waves in the fluid by the velocity potential 
$\phi(x_1, x_2, x_3, t)$ and the pressure $p(x_1, x_2, x_3, t)$. Linearized equations connecting these values are as follows: 
\begin{equation}
\Delta \phi = W ,
\label{e:der1}
\end{equation}
\begin{equation}
\frac{\ptl \phi  }{\ptl t}  = - \frac{1}{\rho}  p, 
\label{e:der2}
\end{equation}
where $\rho$ is the density of the fluid, and 
$W$ is the source function:
\[
W(x_1, x_2, x_3, t) = A v (\delta (x_1-l/2 - vt) - \delta (x_1+l/2 - vt)) \delta(x_2) \delta(x_3) \Theta(t) 
\]
\begin{equation}
\qquad \qquad \qquad \qquad \qquad \qquad
\approx 
- V v \delta'(x_1 - vt) \delta(x_2) \delta(x_3) \Theta(t),
\label{e:der3}
\end{equation}
$\delta$, $\delta'$, and $\Theta$ are the Dirac function, its derivative,  and the Heaviside function, respectively. 

The free boundary condition is
\begin{equation}
\left.
\left( a \rho \frac{\ptl \phi}{\ptl x_3}  - \frac{\ptl p}{\ptl t}  \right) \right|_{x_3 = +0} = 0.
\label{e:der4}
\end{equation}
Here $x_3 = +0$ means that the value is taken above the body. 

Relations (\ref{e:der1}), (\ref{e:der2}), (\ref{e:der4}) can be combined as follows: 
\begin{equation}
\Delta \phi = 0, \qquad x_3 < -0, 
\label{e:der5}
\end{equation}
\begin{equation}
\left. \left( 
\frac{\ptl \phi}{\ptl x_3} + \frac{1}{a} \frac{\ptl^2 \phi}{\ptl t^2} 
\right) \right|_{x_3 = -0} = -  Vv \delta' (x_1 - vt) \delta(x_2) \Theta (t). 
\label{e:der6}
\end{equation}
Here $x_3=-0$ is the position under the body.  

Apply the Fourier--Laplace transform to the solution: 
\begin{equation}
\hat \phi(k_1, k_2, \omega; x_3) = \frac{1}{8 \pi^3}
\int_0^{\infty}
\int \!\!\! \int_{-\infty}^{\infty} 
\phi(x_1,x_2,x_3,t) e^{-i(k_1 x_1 + k_2 x_2 - \omega t)}
dx_1 \, dx_2 \, dt.
\label{e:der7}
\end{equation}
As a result, get an ordinary differential equation 
\begin{equation}
\frac{\ptl^2 \hat \phi}{\ptl x_3^2} = (k_1^2 + k_2^2) \hat \phi, 
\qquad x_3 < -0
\label{e:der8}
\end{equation}
with the boundary condition
\begin{equation}
\left( 
\frac{\ptl}{\ptl x_3} - \frac{\omega^2}{a}
\right)\hat \phi(k_1, k_2, \omega; -0) = -\frac{Vv}{8 \pi^3} \frac{k_1}{\omega - k_1 v}.
\label{e:der9}
\end{equation}
Picking up the solution of (\ref{e:der8}) decaying as $x_3 \to - \infty$, obtain 
\begin{equation}
\phi(k_1, k_2, \omega; -0) = \frac{Vva}{8 \pi^3}
\frac{k_1}{(\omega - k_1 v)(\omega^2 - a \sqrt{k_1^2 + k_2^2})}.
\label{e:der10}
\end{equation}
Finally, we note that
\[
H(x_1, x_2, t) = \frac{p(x_1, x_2, 0, t)}{\rho a} = - \frac{1}{a} \frac{\ptl \phi}{\ptl t} (x_1, x_2, 0, t),
\]
that $\phi$ is continuous at $x_3 = 0$, and apply an inversion of the Fourier transform (\ref{e:der7}). 
The result is (\ref{e:101}).

We note that (\ref{e:101}) can be derived directly from \cite{Liu2001}.

{\bf Preparatory steps.}
Let us evaluate the integral (\ref{e:101}) asymptotically at large values of arguments. For this, we formally write
\[
(x_1, x_2 , t) = (\Lambda X_1 , \Lambda X_2 , \Lambda T),
\]
($\Lambda$ is dimensionless), fix $(X_1, X_2, T)$ and consider the limit $\Lambda \to \infty$. 

Let us introduce dimensionless variables
\[
\xi_1 = \frac{v^2}{a} k_1, \quad \xi_2 = \frac{v^2}{a} k_2, \quad \varpi_3 = \frac{v}{a} \omega,
\quad
z_1 = \frac{a}{v^2} X_1, \quad z_2 = \frac{a}{v^2} X_2, \quad \tau = \frac{a}{v} T. 
\]
In new variables the integral (\ref{e:101}) is rewritten as 
\begin{equation}
H(\Lambda z_1, \Lambda z_2, \Lambda \tau ) = \frac{V a^2 i}{8 \pi^3 v^4} 
 \int \limits_{-\infty + i \epsilon}^{\infty + i \epsilon}
\int \limits_{-\infty}^{\infty }
\int \limits_{-\infty}^{\infty }
\frac{\xi_1 \varpi e^{i \Lambda(\xi_1 z_1 + \xi_2 z_2 - \varpi \tau)}}{(\varpi - \xi_1)(\varpi^2 - \sqrt{\xi_1^2 + \xi_2^2})}
d\xi_1 \, d\xi_2  \, d\varpi.
\label{e:101a}
\end{equation}

The integral (\ref{e:101a}) is well studied. The classical approach \cite{Stoker1957} is to evaluate the integral as a set of nested integrals by performing integration of each of the variables sequentially. Here, however, we would like to demonstrate 
the benefits of our approach that considers the integral as a whole without reducing its dimension.  
As far as we know, the first observation that the global consideration of singularities of multivariate Fourier integrals 
by taking into account the geometry of singularities of the integrand leads to physically clear results was made in \cite{Chapman2002}.

It is easy to see that the integral (\ref{e:101a}) is similar in form to (\ref{e:001}). 
The variables $(\xi_1 , \xi_2 , \varpi)$ correspond to $\xi$, 
and the variables $(z_1 , z_2 , \tau)$ correspond to $z$ (we have changed the notations $\xi_3$ and $z_3$ to make the temporal variables
visually different from the spatial ones).

The plan of the study is as follows: 

\begin{itemize}
\item
The integral is {\em complexified}, i.e.\ 

\item
The integrand is decomposed into multipliers $e^{i \Lambda G}$ and~$F$. 

\item
The singularities $\sigma_j$ of the multiplier~$F$ are found. 

\item
The traces $\sigma'_j$ of these singularities on~$\mathbb{R}^3$ are built. 

\item
The special points are found in $\mathbb{R}^3$  according to our classification.

\item
For special points, additional conditions are checked guaranteeing that the points contribute to the asymptotics.

\item
If a special point gives a contribution to the asymptotics, the asymptotics is calculated using some ready formulas 
or {\em ad hoc}.  

\item 
The positions of the special points $\xi^*$ depend on~$z$.
The phase multiplier $\exp \{ i \Lambda G(\xi^*(z) ; z)\}$ 
is analyzed
for each family of singular points to determine to which type of waves this family corresponds.  

\end{itemize}

Since the integral is nested, 
\[
d\xi_1 \, d\xi_2 \, d\varpi = d\xi_1 \wedge d\xi_2 \wedge d\varpi. 
\]
The functions $F$ and $G$ are written as
\[
G(\xi ; z) = \xi_1 z_1 + \xi_2 z_2 - \varpi \tau,   
\]
\[
F(\xi) =  \frac{i a^2 V}{8 \pi^3 v^4} \frac{\xi_1  \varpi}{(\varpi - \xi_1)(\varpi^2 - \sqrt{\xi_1^2 + \xi_2^2})}
\]

Note that the integral (\ref{e:101}) does not fully satisfy the constraints imposed above. Near the origin, the function $F$ has branch sets defined by equations
\[
\xi_1 \pm i \xi_2 = 0,
\]
which have no real property. These sets intersect the integration domain (nevertheless, the integral (\ref{e:101}) is defined correctly). This circumstance does not allow for free deformation of the integration domain near the origin. Such points require additional theoretical work, which is beyond the scope of this paper, i.e., we will not calculate here the asympothetic terms generated by the neighborhood of origin in the space~$\xi$. 

\begin{figure}[ht]
\centerline{\epsfig{width = 8 cm, file=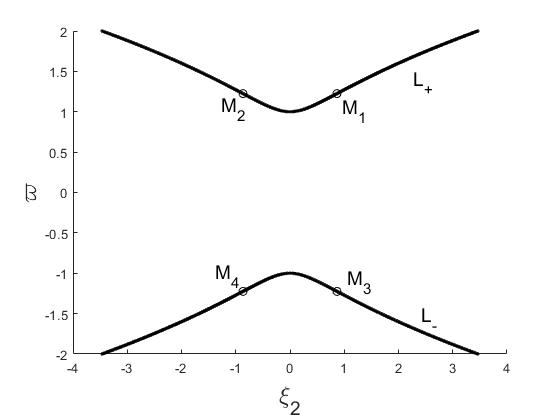}}
\caption{Curves $L_\pm$ }
\label{fig04}
\end{figure}

\vskip 6pt
\noindent
{\bf Geometry of singularities.} 
Let us apply the technique to evaluate the integral developed above. The function $F$, considered in $\mathbb{R}^3_\delta$ without a neighborhood of the origin, has two polar singularities $\sigma_1$ and $\sigma_2$, given by the functions
\begin{equation}
g_1(\xi_1, \xi_2 , \varpi) = \varpi - \xi_1, 
\label{e:102}
\end{equation}
\begin{equation}
g_2(\xi_1, \xi_2 , \varpi) = \varpi^2  - \sqrt{\xi_1^2 + \xi_2^2}.  
\label{e:103}
\end{equation}
Both singularities are bypassed from above with respect to the variable~$\varpi$.
The surfaces $\sigma'_1$ and $\sigma_2'$ are a plane and a curvilinear cone. 
There are two lines of intersection of these surfaces: $L = L_1 \cup L_2 = \sigma'_1 \cap \sigma'_2$, 
parametrized by the variable $\varpi$ as
\begin{equation}
L = 
\{ (\xi_1 , \xi_2, \varpi) \in \mathbb{R}^3
\, \, : \, \, 
\xi_1 = \varpi , \quad
\xi_2 = \pm \sqrt{\varpi^4 -\varpi^2}
\}, 
\qquad 
|\varpi| \ge 1.
\label{e:104}
\end{equation}
The shape of these curves is shown in \figurename~\ref{fig04} 
and it 
is important for the consideration below. 
These curves (as curves in the plane $\sigma'_1$) have inflection points $M_1 , \dots , M_4$,
corresponding to $\varpi = \pm \sqrt{3/2}$.

According to our classification, the following special points can be expected to appear: stationary phase points on the singularities
$\sigma'_1$ and $\sigma'_2$, and stationary phase points on the crossing 
$\sigma_1' \cap \sigma_2'$. Furthermore, as we have noted, the neighborhood of the origin does not fall within our classification, and it needs to be considered separately. 

Note that there cannot be a stationary phase point in 3D since $G$ is a linear function, and there are no triple crossings of singularities since there are only two singularities. 

\vskip 6pt
\noindent
{\bf Stationary phase points at $\sigma_1' \cap \sigma_2'$. Kelvin waves. } Here we show that  
the most interesting wave in the system, namely the wedge-shaped wave trace of a moving body, is given by stationary phase points at the
intersection of singularities. Moreover, the basic properties of Kelvin waves can be deduced from geometric considerations. 

Consider the pair of curves $L$ introduced above. Each point of $L$ at certain coordinates $z$ can be a stationary phase point at $\sigma'_1 \cap \sigma'_2$. Let us find these coordinates~$z$. Let $\xi^* \equiv (\xi_1^*, \xi_2^*, \varpi_*) \in L$. 
Differentiating (\ref{e:104}) by $\varpi$,
find the vector ${\bf a}$ tangent to $L$ at the point~$\xi^*$:
\begin{equation}
{\bf a} = \left( 
1, \frac{2 \varpi_*^2 - 1}{\sqrt{\varpi_*^2 - 1}} , 1
\right) ,
\label{e:105}
\end{equation}
where the square root takes one of two values depending on the sign in (\ref{e:104}). 
Condition (\ref{e:016}) is written as 
\begin{equation}
z \cdot {\bf a} = z_1 + z_2 \frac{2 \varpi_*^2 - 1}{\sqrt{\varpi_*^2 - 1}} - \tau = 0 .
\label{e:107}
\end{equation}

Rewrite (\ref{e:107}) as 
\begin{equation}
\frac{z_2}{\tau - z_1} = \frac{\sqrt{\varpi^2 -1}}{2 \varpi^2 -1},
\label{e:107b}
\end{equation}
introduce a combination of coordinates 
\begin{equation}
\lambda = \frac{z_2}{\tau - z_1},
\label{e:107a}
\end{equation}
and solve equation (\ref{e:107}) with respect to~$\varpi$:
\begin{equation}
\varpi_*(\lambda) = \pm \frac{\sqrt{4 \lambda^2 + 1 \pm \sqrt{1 - 8 \lambda^2}}}{2\sqrt{2} \lambda} , 
\qquad 
\label{e:108}
\end{equation}
Substituting (\ref{e:108}) into (\ref{e:104}) gives a point on~$L$.

Let us analyze formula (\ref{e:108}). First, an obvious condition for the 
existence of a real solution is $1 - 8 \lambda^2 \ge 0$.  This condition can be rewritten as
\begin{equation}
\left| \frac{z_2}{\tau - z_1}\right| \le  \frac{1}{2 \sqrt{2}}.
\label{e:109}
\end{equation}
Let us assume that $\tau - z_1 >0$ (we will comment on this condition below). Then
(\ref{e:109}) defines the angular region behind the current position ($z_1 = \tau$) of the source. Note that the angle 
\begin{equation}
\theta = {\rm arctan} \left(1/(2 \sqrt2)\right) \approx 19.47^\circ
\label{e:109a}
\end{equation}
is the famous Kelvin angle containing waves following a ship. 

\begin{figure}[ht]
\centerline{\epsfig{width = 12 cm, file=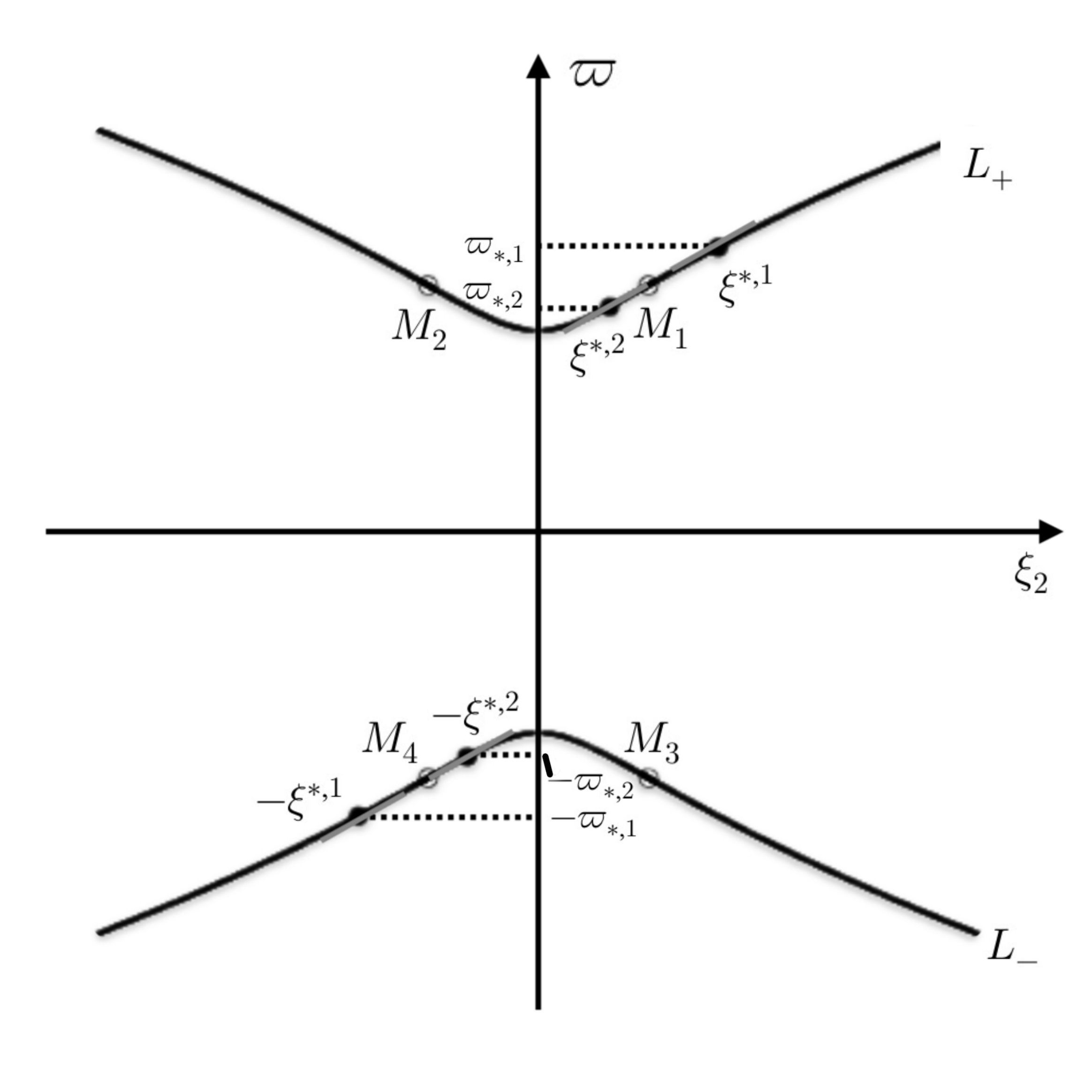}}
\caption{Stationary phase points at $L$ and the wave behind the body.}
\label{fig05}
\end{figure}

If condition (\ref{e:109}) is fulfilled  strictly, the formula gives four points
$\xi^*$ on $L$ denoted as $\xi^{*,1}$, $\xi^{*,2}$, $-\xi^{*,1}$, $-\xi^{*,2}$.
For these points 
$\varpi_* = \pm \varpi_{*,1}, 
\pm \varpi_{*,2}$,
\begin{equation}
\varpi_{*,1}(\lambda) = \frac{\sqrt{4 \lambda^2 + 1  +  \sqrt{1 - 8 \lambda^2}}}{2\sqrt{2} \lambda} , 
\quad 
\varpi_{*,2}(\lambda) = \frac{\sqrt{4 \lambda^2 + 1  -  \sqrt{1 - 8 \lambda^2}}}{2\sqrt{2} \lambda} .
\label{e:109b}
\end{equation}
At $z_2 > 0$ and $\tau - z_1 > 0$ the position of the points is shown in \figurename~\ref{fig05}. Obviously, the tangents to $L$ at all four points are parallel.

In the following, we will only consider contributions from points $\xi^{*,j}$, $j = 1,2$. An elementary analysis shows that the contributions of the points $- \xi^{*,j}$ are complex conjugate to the contributions from~$\xi^{*,j}$.  

\begin{figure}[ht]
\centerline{\epsfig{width = 7 cm, file=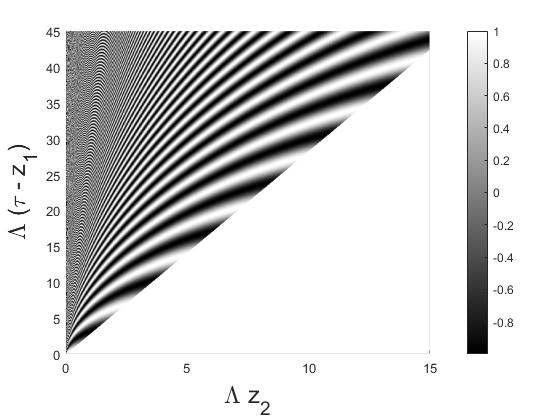}\epsfig{width = 7 cm, file=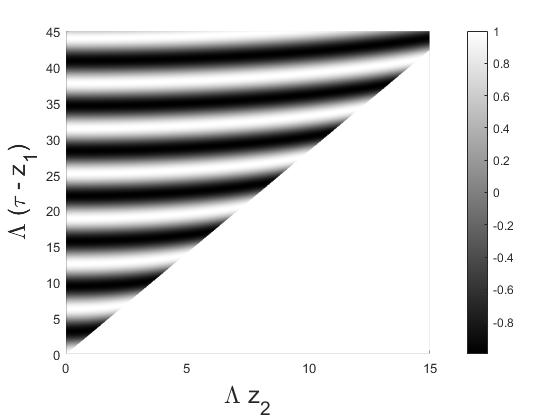}}
\caption{The function $\cos (\Lambda G(\xi^{*,1}))$ (left) и $\cos (\Lambda G(\xi^{*,2}))$ (right)}
\label{fig06}
\end{figure}

The point $\xi^{*,1}$ corresponds to one type of the waves, and $\xi^{*,2}$ corresponds to another one. 
Let us plot a picture of wave fronts of each type. According to (\ref{e:017}), the wave fronts are defined by the factor 
$\exp\{ i \Lambda G(\xi^*)\}$ (all other multipliers give only information about the amplitude and initial phase). Take
\begin{equation}
G(\xi^{*,j}) =  \frac{\varpi_{*,j}^3(\lambda) (z_1 - \tau)}{2 \varpi_{*,j}^2 (\lambda) -1},
\label{e:111}
\end{equation}
and plot the graph of the real part of $\exp \{i \Lambda G(\xi^*)\}$ 
(see \figurename~\ref{fig06}) to get a picture of the wave fronts. It is easy to see that the exponential factor $\exp\{ i \Lambda G(\xi^*) \}$
correctly describes the phase in the two
families of Kelvin waves that follow the body.

At $8 \lambda^2 = 1$ points $\xi^{*,1}$ and $\xi^{*,2}$ merge at the inflection points of the curve~$L_+$.
The merging of the points $\xi^{*,1}$ and $\xi^{*,2}$ corresponds to the boundary of the area occupied by the wave fronts. 

Let us find the condition that the point $\xi^*$ contributes to the asymptotics of the integral. For this, we need to check that $\Gamma$ passes below the singularities $w_1 = 0$ and $w_2 = 0$. 
It follows from (\ref{e:101a}) that the
domain $\Gamma$ is shifted upward with respect to the variable~$\varpi$. 
Hence, it bypasses above the set $\varpi - \xi = 0$.
The condition of bypassing  below $w_1 = 0$ is written
as $\alpha_1 < 0$ or 
\[
z_1 > \frac{\tau - z_1}{2 \varpi_{*,j}^2 -1} .
\]
Physically, this condition is related to the fact that the trace starts to form only at the moment $\tau = 0$,  and if the condition 
is not fulfilled, the corresponding part of the trace does not have time to form. 

The domain $\Gamma$ bypasses above the set $\varpi^2 - \sqrt{\xi_1^2 + \xi^2_2}$ for $\varpi > 0$ and below this set for $\varpi < 0$. Thus, 
the condition for $\Gamma$ bypassing below $w_2 = 0$ 
is equivalent to the fact that $\alpha_2 \varpi_{*,j} < 0$ or, according to (\ref{e:107}),
\[
\tau - z_1 > 0.
\]
This is easy to interpret physically: the trail should be behind the body. 

Thus, a purely geometric consideration within the framework of the proposed method allowed us to establish that two types of waves exist in 
a certain part of the $z$-space and to construct the corresponding wave picture.

Let us find the wave amplitudes amplitude information using (\ref{e:017}).
Consider the integral of the region near the point $\xi^{*,j}$ in more detail. We will assume that this point is not close to the inflection point $L$ (the neighborhood of the inflection point should be considered separately). In accordance with (\ref{e:016a}), (\ref{e:016b}), introduce the coordinates
\begin{equation}
w_1 = \alpha_1 (\varpi - \xi_1),
\qquad 
w_2 = \alpha_2 \left( \varpi^2 - \sqrt{\xi_1^2 + \xi_2^2} \right),
\qquad 
w_3 = \varpi - \varpi_{*,j}(\lambda). 
\end{equation}
By selecting
\begin{equation}
\alpha_1 = -z_1 + \frac{\tau - z_1}{2 \varpi^2_{*,j} (\lambda) -1} , 
\qquad 
\alpha_2 = - \frac{z_2 \varpi_{*,j}}{\sqrt{\varpi^2_{*,j} - 1}}
\label{e:110}
\end{equation}

function $G$ is written as 
\begin{equation}
G \approx G(\xi^{*,j}) + w_1 + w_2 + \frac{1}{2} \beta w_3^2,
\label{e:110a}
\end{equation}
 
\begin{equation}
\beta = \frac{\varpi_{*,j} (\lambda) (2  \varpi^2_{*,j} (\lambda) - 3)}{\left(\varpi_{*,j}^2 (\lambda) -1\right)^{3/2}}.
\label{e:112}
\end{equation}
The amplitude multiplier in (\ref{e:107}) can be derived using (\ref{e:019a}) with $\mu_1 = \mu_2 = -1$. 
The determinant of the Jacobian is easy to compute: 
\[
J=\frac{1}{z_2}\left(z_1-\frac{\tau-\tau_1}{2 \varpi^2_{*,j} (\lambda) -1}\right)^{-1}.
\]
The final expression is rather lengthy, and we do not put it here. 

\vskip 6pt
\noindent
{\bf Stationary phase point on the surface $\sigma'_2$. The transient process caused by the onset of body motion.}
Let us choose the point $\xi^* = (\xi^*_1, \xi_2^*, \varpi_*) \in \sigma'_2$, defined by the real parameters $\varpi_* , \varphi$ as
\begin{equation}
\xi^*_1 = \varpi_*^2 \cos \varphi, 
\qquad 
\xi^*_2 = \varpi_*^2 \sin \varphi. 
\label{e:113a}
\end{equation}
The normal vector ${\bf b}$ to $\sigma_2'$ at this point is
\[
{\bf b} = (- \cos \varphi , - \sin \varphi , 2 \varpi_*).
\]
Condition (\ref{e:015}) means that the vectors $\nabla G = (z_1, z_2, -\tau)$ and ${\bf b}$ are collinear, i.e.\ 
\begin{equation}
(z_1, z_2)  = \frac{\tau}{2 \varpi_*} (\cos \varphi,  \sin \varphi). 
\label{e:113}
\end{equation}
Obviously, $\sigma_2'$ is the dispersion curve for gravitational waves
in a deep fluid. The group velocity at frequency $\varpi_*$ is $1/(2 \varpi_*)$, 
therefore (\ref{e:113}) describes the set of points that the pulse 
having frequency $\varpi_*$
reaches in time~$\tau$.
Note that the point with coordinates $(z_1, z_2, \tau)$ corresponds to a frequency
\begin{equation}
\varpi_* (z_1 , z_2 , \tau ) = \frac{\tau}{2 \sqrt{z_1^2 + z_2^2}}
\label{e:113b}
\end{equation}
Thus, for almost any point $z = (z_1, z_2, \tau)$ there is a stationary point $\xi^*$ on $\sigma'_2$, and formulas 
(\ref{e:113b}), (\ref{e:113}), and (\ref{e:113a}) allow us to find this point.  

According to (\ref{e:015z}), we introduce the coordinates 
\[
w_1 = \alpha \left(\varpi^2 - \sqrt{\xi_1^2 + \xi_2^2} \right) ,
\]
\[ 
w_2 = (\xi_1 - \varpi_*^2 \cos\varphi) \cos \varphi + (\xi_2 - \varpi_*^2 \sin\varphi) \sin \varphi , 
\]
\[
w_3 = (\xi_1 - \varpi_*^2 \cos\varphi) \sin \varphi - (\xi_2 - \varpi_*^2 \sin\varphi) \cos \varphi . 
\]
For
\begin{equation}
\alpha = - \frac{\tau}{2 \varpi_*} = - \sqrt{z_1^2 + z_2^2},
\label{e:114}
\end{equation}
the phase function $G$ can be written as 
\begin{equation}
G \approx G(\xi^*) + w_1 + \frac{1}{2} (\beta_2 w_2^2 + \beta_3 w_3^2) , 
\label{e:115}
\end{equation}
where 
\begin{equation}
G(\xi^*) = - \frac{\tau^2}{4 \sqrt{z_1^2 + z_2^2}} ,
\label{e:116}
\end{equation}
\begin{equation}
\beta_2 = \frac{\tau}{4 \varpi_*^3} = \frac{2(z_1^2 + z_2^2)^{3/2}}{ \tau^2}, 
\qquad 
\beta_3 = -\frac{\tau}{2 \varpi_*^3} = - \frac{4(z_1^2 + z_2^2)^{3/2}}{ \tau^2}
\label{e:118}
\end{equation}

The condition that the point $\xi^*$ contributes to the 
asymptotics is that the original integration domain passes below $w_1 = 0$. 
As before, this condition is given by the inequality 
$\alpha \varpi_* < 0$ that is equivalent to the inequality $\tau > 0$ guaranteeing causality.

Let us rewrite $F$ in the form (\ref{e:015y}) near $\xi^*$: 
\[
F(w_1, w_2, w_3) \approx C w_1^{-1},
\]
\begin{equation}
C 
= \frac{Va^2 i }{16 \pi^3 v^4}\frac{\tau^2 \cos \varphi}{\tau \cos \varphi - 2 \sqrt{z_1^2 + z_2^2} }.
\label{e:117}
\end{equation}

The asymptotics of the corresponding wavefield component is given by substituting (\ref{e:117}), (\ref{e:118}), (\ref{e:116})  into (\ref{e:015a}). The determinant of the Jacobian $J$ is easily calculated:
\[
J = - \tau^{-1}.
\]

We conclude that the contribution of the stationary phase points on $\sigma'_2$ is a cylindrical wave going from the origin. 
The source of this field is a concentrated force acting at the
moment $\tau = 0$. We can interpret this contribution as a transient process associated with the onset of body motion. 
Note that the wave field depends on the angle~$\varphi$, i.e.\ there is no rotational symmetry. 

The denominator (\ref{e:117}) is singular for some~$z$. 
This means that the asymptotic formula (\ref{e:015a}) does not work. The reason is that for such $z$ the stationary phase point on $\sigma_2'$ merges with the stationary phase point on $\sigma'_1 \cap \sigma'_2$. In this case, the following
conditions of using both asympotics are violated, and  it is necessary to construct a 
more complicated transient asymptotics in the neighborhood of such points.

\vskip 6pt
\noindent
{\bf The remaining special points.}
In this paper, we omit the consideration of a number of special points. 
That is, at
$z_2 = 0$, $z_1 = \tau$ all the points of the plane $\sigma'_1$ become points of the stationary phase on~$\sigma'_1$. 
Thus, the locality  principle does not work, and it is necessary to calculate a two-dimensional
integral over the whole~$\sigma'_1$. This integral will describe the field accompanying the moving body. 

Moreover, we do not consider the neighborhood of origin in the space $\xi$. As we mentioned above, this point does not fall under the special point classification made in the paper. However,  we can construct the asymptotic contribution using similar methods for this special point. The corresponding contribution gives the field behavior at low frequencies. 

Finally, we do not consider the process of merging several special points (most interesting would be the merging of $\xi^{*,1}$ и $\xi^{*,2}$). This consideration requires leaving more power terms in the Taylor series for $G$. Obviously, the
transition zone will be described by the Airy function.

\section{Conclusion}

The authors propose a technique for constructing asymptotic estimates of three-dimensional Fourier and Fourier-type 
integrals for analytic integrands. 
The main attention is paid to the influence of singularities of the integrand. It is shown
that in the general case the principle of locality works, i.e., the evaluation of the integral is determined by several special points, the main types of which are listed in the paper. It is important to note that special points and points where the integrand is singular are different sets.

For some types of singular points, there are geometric conditions under which a given point contributes to the integral. These conditions reflect the mutual position of the initial integration domain and the singularities.

The leading terms of the asymptotic expansions of the integral are constructed in general form for the main types of special points. Whenever a convenient local coordinate system is constructed, the integral is reduced to a nested integral in which, for some cases, the integral is a saddle-point integral and for the rest it is of the Gamma-function type.

The described technique is tested on the classical problem of Kelvin waves on the surface of a deep liquid. The waves behind the 
moving body are shown to be described by one of the types of special points: the stationary phase points at the intersection of singularities. The known angle (\ref{e:109a}) at which these waves are concentrated is related to the position of 
the inflection point on the line of intersection of singularities. Transient processes are taken into account easily using the proposed method. 

\section*{Acknowledgements}

Authors are grateful to T.~Shugailo, Prof.\ G.~Mishuris, Dr.~A.~Korolkov, and Dr.~R.~Assier for fruitful discussions.  

The research was supported by Russian Science Foundation grant 25--22--00106, https://rscf.ru/project/25-22-00106/.


\begin{thebibliography}{99}

\bibitem{Martin2025} 
P. A. Martin, Kelvin's method of stationary phase? \textit{Wave Motion}, 134:103481, 2025. 


\bibitem{Fedoryuk1977} 
M. V. Fedoryuk, \textit{The Method of Steepest Descent}, Nauka, Moscow, 1977 (in Russian).

\bibitem{Mironov2021}
M. A. Mironov, A. V. Shanin, A. I. Korolkov, K. S. Kniazeva, Transient processes in a gas/plate structure in the case of light loading. \textit{Proceedings of the Royal Society A: Mathematical, Physical and Engineering Sciences}, 477:20210530, 2021.

\bibitem{Assier2022}
R. C. Assier, A. V. Shanin, A. I. Korolkov, A contribution to the mathematical theory of diffraction: A note on double fourier integrals. \textit{Quarterly Journal of Mechanics and Applied Mathematics}, 76(1):1-47, 2022. 

\bibitem{Shanin2024}
A. V. Shanin, R. C. Assier, A. I. Korolkov, O. I. Makarov, Double Floquet-Bloch transforms and the far-field asymptotics of Green's functions tailored to periodic structures. \textit{Physical Review B}, 110(2):024310, 2024.

\bibitem{Poincare1904}
H. Poincare, Sur la diffraction des ondes electriques: à propos d'un article de M. Macdonald. \textit{Proceedings of the Royal Society of London}, 72(477-486):42-52, 1904.

\bibitem{Thomson1887}
W. Thomson, On ship waves. \textit{Proceedings of the institution of mechanical engineers}, 38(1):409-434, 1887.

\bibitem{Kelvin1906}
Lord Kelvin, Deep sea ship-waves. \textit{Proceedings of the royal society of Edinburgh}, 25(2):1060-1084, 1906.

\bibitem{Lamb1916a}
H. Lamb, LXV. On wave-patterns due to a travelling disturbance. \textit{The London, Edinburgh, and Dublin Philosophical Magazine and Journal of Science}, 31(186):539-548, 1916.

\bibitem{Lamb1916b}
H. Lamb, \textit{Hydrodynamics},  New York, Dover publications, 1945.

\bibitem{Stoker1957}
J. J. Stoker, \textit{Water waves: The mathematical theory with applications}, Interscience Publishers, Inc., New York, 1957.

\bibitem{Wehausen1960}
J. V. Wehausen, E. V. Laitone, \textit{Surface waves}, Fluid Dynamics/Strömungsmechanik, Berlin, Heidelberg: Springer Berlin Heidelberg, 1960. P. 446-778.

\bibitem{Liu2001}
M. Liu, M. Tao, Transient ship waves on an incompressible fluid of infinite depth. \textit{Physics of Fluids}, 13(12):3610-3623, 2001.

\bibitem{Shabat2}
B.V. Shabat. \textit{Introduction to complex analysis Part II. Functions of several variables}. AMS, 1992.

\bibitem{Huybrechs2007}
D. Huybrechs, S. Vandewalle, The construction of cubature rules for multivariate highly oscillatory integrals. \textit{Mathematics of computation}, 76(260):1955-1980, 2007.



\bibitem{Chapman2002}
C. J. Chapman,  The wavenumber surface in blade--wortex interaction. \textit{Proceedings of the IUTAM symposium on diffraction and scattering in fluid mechanics and elasticity, Manchester, UK 2000, ed.\ by I.D.Abrahams, P.A.Martin, M.J.Simon, Kluwer}, 169-178, 2002. 


\end{thebibliography}
\end{document}